# Two Fefferman-Type Constructions Involving Almost Grassmann Structures and Path Geometries

**Zhangwen Guo**[1]

Received: 8 October 2025 / Accepted: 25 September 2026

## Abstract

We introduce a Fefferman-type construction that associates an almost Grassmann structure of type $(2, n+1)$ to every $(n+1)$-dimensional path geometry. We prove that the construction is normal and provide two equivalent characterizing conditions for all almost Grassmann structures which locally arise from this construction: one in terms of certain parallel tractors and the other in terms of a Weyl connection of an almost Grassmann structure. We prove that the latter condition is independent of the choice of Weyl connection. We then introduce a related Fefferman-type construction associating an almost Grassmann structure of type $(2, n+1)$ to every almost Grassmann structure of type $(2, n)$. We prove that this construction is non-normal and characterize all almost Grassmann structures which locally arise in this way in Cartan geometric terms.

## 1 Introduction

Assume that $N^{n+1}$ is an $(n+1)$-dimensional manifold, $n \geq 1$. A path geometry on $N$ describes a smooth collection of paths, i.e., immersed 1-dimensional submanifolds, on $N$ such that there is a unique such path for every point and associated direction. On a formal level, the paths are encoded on the projectivized tangent bundle $M^{2n+1} = \mathcal{P}TN$ as the unique line distribution $E \subseteq TM$ such that (i) the rank $(n+1)$-tautological

The author is grateful to her doctoral advisor A. Čap for introducing her to the problems solved in this paper and his steady support. This research was funded in part by the Austrian Science Fund (FWF): 10.55776/P33559 and 10.55776/Y963.

✉ Zhangwen Guo
zhangwen.guo@univie.ac.at

1 Faculty of Mathematics, University of Vienna, Oskar-Morgenstern-Platz 1, 1090 Vienna, Austria

distribution of $M$ can be written as the direct sum of $E$ and the vertical bundle $V \subseteq TM$ of $M \to N$ and (ii) the integral manifolds of $E$ descend to the paths on $N$. In particular, there is a construction $(N, \{\text{paths}\}) \mapsto (M, E, V)$ leading to a $(2n+1)$-dimensional generalized path geometry; see Sect. 4.2. Conversely, it has been proven that if the defining distribution $V$ of a $(2n+1)$-dimensional generalized path geometry is involutive, then this generalized path geometry locally arises from the construction described above. For this reason, we will refer to a generalized path geometry with involutive $V$ as a path geometry. Moreover, the defining rank-$n$-distribution has been shown to always be involutive if $n \neq 2$; see, e.g., [11, Proposition 4.4.4]. In particular, a generalized path geometry $(M, E, V)$ with $\dim(M) \neq 5$ is always a path geometry.

Assume now that $M^{pq}$ is a $pq$-dimensional manifold, $p, q \geq 2$. We follow the convention that an almost Grassmann structure $(M, E, F)$ of type $(p, q)$ on $M$ is defined by two auxiliary vector bundles $E \to M$ and $F \to M$ of rank $p$ and of rank $q$, respectively, and isomorphisms

$$E^* \otimes F \cong TM \quad \text{and} \quad \wedge^p E^* \cong \wedge^q F. \tag{1}$$

Recall from, e.g., [11, Subsections 4.1.3, 4.2.3 and 4.4.3] the following background. Each of the above geometric structures carries a canonical parabolic geometry, which in turn comes with its invariants. The types of fundamental invariants splits into different cases. For the general case of a $(2n+1)$-dimensional generalized path geometry with either $n = 2$ and $V$ involutive or $n \geq 3$ (where $V$ is automatically involutive), the fundamental invariants consist of a harmonic torsion and a harmonic curvature; see (18) below. In the case where $n = 2$ and $V$ not involutive, there is an extra invariant torsion in addition to the general case. For $n = 1$, the fundamental invariants consists of two Cotten-York-type tensors; see Sect. 5.1 below. On the other hand, for the general case of an almost Grassmann structure of type $(2, n)$, $n \geq 3$, the fundamental invariants consist of a harmonic torsion and a harmonic curvature; see (11) below. The fundamental invariants for an almost Grassmann structure of type $(2, 2)$ is a pair of harmonic curvatures; see Sect. 5.1. Finally, the fundamental invariants for an almost Grassmann structure of type $(p, q)$, $p, q \geq 3$ is a pair of harmonic torsions.

In this present article, we only deal with the general situation of the two geometric structures described above, with the exception of Sect. 5.1 where 3-dimensional generalized path geometries and almost Grassmann structures of type $(2, 2)$ are discussed. In particular, the harmonic torsion $\tau \in \Omega^2(M, TM)$ of an almost Grassmann structure $(M, E, F)$ of type $(2, n)$, $n > 2$, coincides with the torsion of any Weyl connection on $(M, E, F)$. Indeed, the Weyl connections on $(M, E, F)$ are exactly pairs of connections on $(E, F)$ which are compatible with the second isomorphism in (1) and, via the first isomorphism in (1), induce a connection on $TM$ whose torsion is $\partial^*$-closed and hence equal to $\tau$; see, e.g., [11, Theorem 5.2.3 (1)]. Since $\tau$ is $\partial^*$-closed, it vanishes upon any contraction. Hence, up to a sign, there is a unique nontrivial contraction of $\tau \otimes \tau$ resulting in a $(0, 2)$-tensor field. This contraction is clearly symmetric. We will denote it by

$$tr(\iota_\tau \tau) \in \Gamma(\mathrm{Sym}^2\, T^*M). \tag{2}$$

We introduce a Fefferman-type functor that associates to every path geometry $(M^{2n+1}, E, V)$ an almost Grassmann structure of type $(2, n+1)$. This is motivated by the work of Crampin and Saunders [14], who associate to every path geometry $(M^{2n+1}, E, V)$ a Segre structure of type $(2, n+1)$. Recall that, by definition, every almost Grassmann structure induces a Segre structure of the same type; see, e.g., [25].

The framework of Fefferman-type functors is developed from the pioneering work of Fefferman [16], who associates to every strictly pseudoconvex real hypersurface in $\mathbb{C}^n$ a conformal geometry on the total space of a circle bundle. It has been generalized to a construction associating to every abstract CR geometry a conformal geometry by Burns, Diederich and Shnider in [3]. As a generalization of the latter, the Fefferman-type functors on Cartan geometries have been introduced by Čap in [5]. The framework of Fefferman-type functors has been extended to a slightly broader class of functors with the same properties as the original one by Čap and Žádník in [13]. This is based on the more general framework of so-called extension functors on Cartan geometries.

The Fefferman-type functor we specify is applied to the canonical parabolic geometry associated to every path geometry $(M^{2n+1}, E, V)$ together with an auxiliary orientation on $E \oplus V$ or, equivalently, on $E^* \otimes \wedge^n V$, since $\otimes^2 E$ is canonically oriented. The result is a parabolic geometry that induces an almost Grassmann structure of type $(2, n+1)$ on the principal $\mathbb{R}_+$-bundle $GL_+(E^* \otimes \wedge^n V)$. If the initial path geometry arises from a collection of paths on $N^{n+1}$ as described in the first paragraph in this section, the base space of the resulting almost Grassmann structure is the weighted slit tangent bundle $TN \otimes \mathcal{E}(-1) \setminus \{0\}$. It is precisely the base space of the Segre structure that Crampin and Saunders associate to a path geometry on $N$; see [14, Section 3].

We prove that the resulting parabolic geometry coincides with the parabolic geometry canonically associated to the resulting almost Grassmann structure. In other words, the Fefferman-type construction we specify is normal. In particular, the functor induces an explicit relation between the canonical Cartan connection of a path geometry and the canonical Cartan connection of the resulting almost Grassmann structure. From this relation, it is straightforward to see that the harmonic torsion $\tau$ of every resulting almost Grassmann structure has the property that $tr(\iota_\tau \tau) = 0$. Moreover, the resulting almost Grassmann structure has vanishing harmonic torsion if and only if the initial path geometry has vanishing harmonic torsion.

The relations we find also enable us to characterize all almost Grassmann structures that are locally obtained from the construction. The characterization of conformal geometries that locally arise from the construction of Burns, Diederich and Shnider in [3] has been given by Sparling; see [17, Theorem 3.1]. Sparling's condition characterizing such structures is the existence of a certain isotropic conformal Killing field. In the context of parabolic geometry, such a characterization has been obtained by Čap and Gover in terms of the existence of a parallel tractor; see [7, Theorem 3.1] together with [8, Theorem 2.5]. We also refer the reader to the following articles for some analogous characterizations of geometric structures arising from the respective Fefferman-type constructions.

We also extend our Fefferman-type functor to define a construction from almost Grassmann structures of type $(2, n)$ to almost Grassmann structures of type $(2, n+1)$. We prove that the resulting almost Grassmann structure of type $(2, n+1)$ has the

| Article | Fefferman-type constructions |
|---|---|
| [2] | $(n, n(n+1)/2)$-distribution-to-almost spinorial |
| | (3, 6)-distribution-to-conformal |
| | integrable split signature CR-to-(3, 6)-distribution |
| | integrable Lagrangian contact-to-(3, 6)-distribution |
| [19] | (2, 3, 5)-distribution-to-conformal |
| [1] | quaternionic contact-to-conformal |
| [20] | (2, 3, 5)-distribution-to-conformal |
| | (3, 6)-distribution-to-conformal |
| [23] | (2, 3, 5)-distribution-to-Lie contact |
| [21] | projective-to-conformal |
| | integrable Lagrangian contact-to-conformal |

property $tr(\iota_{\tilde{\tau}}\tilde{\tau}) = 0$ for the harmonic torsion $\tilde{\tau}$ if and only if the initial almost Grassmann structure of type $(2, n)$ has the same property on its harmonic torsion. Moreover, the resulting almost Grassmann structure of type $(2, n+1)$ is torsion-free if and only if the initial almost Grassmann structure of type $(2, n)$ is torsion-free. We note that, among all the above mentioned characterizations of Fefferman-type constructions, that of projective-to-conformal [21] is the only one that deals with so-called non-normal Fefferman-type construction. Our Fefferman-type construction from almost Grassmann structures of type $(2, n)$ to almost Grassmann structures of type $(2, n+1)$ is also of this kind. We follow a procedure similar to that in [21] for characterizing the construction.

The general theory of Fefferman-type constructions prior to the present work has only been well-established on the Cartan level: Given two Cartan geometries related by a Fefferman-type construction, the relation of their Cartan connections and the relation of their Cartan curvatures have been known; moreover, there is an intrinsic characterization in terms of tractor connection and Cartan curvature. These results are reviewed in Sect. 3 and have originated in [5, 9] and [4]. The difficulty of passing these results to a Fefferman-type construction between geometric structures is due to the lack of general theory when the relevant Cartan geometries are parabolic and need to be considered to vary within certain equivalent classes. In particular, essential questions such as normality and normalization have to be investigated on an individual basis. Our technical solution to the above difficulty is to provide case-by-case descriptions on the Kostant codifferentials applied to relevant related tensors; see the proofs of Proposition 15, Theorem 17, Proposition 24 and Theorem 27. Moreover, the proofs of our two characterization results Theorem 17 and Theorem 27 involve normalization in the reduction procedure. In the situation of Theorem 17 the normalization procedure is straightforward, whereas Theorem 27 involves a normalization homogeneity-by-homogeneity, using a general result Proposition 3 originated in [21]. In particular, the current work together with previously mentioned works on Fefferman-type constructions provide ingredients on extending the formal theory to the situation of parabolic geometries.

## 2 A Generalized Correspondence Space Theorem

In this section, we assume that $G$ is a Lie group and $H \subseteq P \subseteq G$ are closed Lie subgroups. Denote their Lie algebras by $\mathfrak{h} \subseteq \mathfrak{p} \subseteq \mathfrak{g}$, respectively.

**Definition 1** (*See, e.g., [11, Subsection 1.5.1]*) A Cartan geometry $(\mathcal{G}, \omega)$ of type $(G, P)$ is a principal $P$-bundle $\mathcal{G} \to M$ together with a Cartan connection $\omega \in \Omega^1(\mathcal{G}, \mathfrak{g})$, i.e., a $\mathfrak{g}$-valued one-form on $\mathcal{G}$ such that (i) $\omega$ is $P$-equivariant with respect to the principal right action and the adjoint action of $G$ on $\mathfrak{g}$ restricted to $P \subseteq G$, (ii) $\omega(\zeta_X) = X$ where $\zeta_X \in \mathfrak{X}(\mathcal{G})$ is the fundamental vector field generated by $X \in \mathfrak{p}$ and (iii) $T\mathcal{G} \xrightarrow[\cong]{(\pi,\omega)} \mathcal{G} \times \mathfrak{g}$ where $\pi : T\mathcal{G} \twoheadrightarrow \mathcal{G}$ is the natural projection.

The curvature of the Cartan connection is the two-form $\kappa \in \Omega^2(\mathcal{G}, \mathfrak{g})$ given by

$$\begin{aligned}\kappa(\xi, \eta) &= d\omega(\xi, \eta) + [\omega(\xi), \omega(\eta)] \\ &= \xi.\omega(\eta) - \eta.\omega(\xi) - \omega([\xi, \eta]) + [\omega(\xi), \omega(\eta)]\end{aligned}$$

for all $\xi, \eta \in \mathfrak{X}(\mathcal{G})$.

It is easy to see that the Cartan connection $\omega : T\mathcal{G} \to \mathfrak{g}$ induces an isomorphism

$$TM \cong \mathcal{G} \times_P \mathfrak{g}/\mathfrak{p}.$$

Recall from, e.g., [11, Lemma 1.5.1] that the curvature $\kappa$ is $P$-equivariant and horizontal. In particular, $\kappa$ may also be expressed as a $P$-equivariant map $\kappa : \mathcal{G} \to L(\wedge^2 \mathfrak{g}/\mathfrak{p}, \mathfrak{g})$ and as a two-form $\kappa \in \Omega^2(M, \mathcal{G} \times_P \mathfrak{g})$.

The orbit space $\tilde{M} = \mathcal{G}/H$ is called the correspondence space of the Cartan geometry $(\mathcal{G} \to M, \omega)$. It canonically carries the Cartan geometry

$$(\mathcal{G} \to \tilde{M}, \omega)$$

of type $(G, H)$, whose curvature is also $\kappa \in \Omega^2(\mathcal{G}, \mathfrak{g})$. In particular, as a section of $\Omega^2(\tilde{M}, \mathcal{G} \times_H \mathfrak{g})$, $\kappa$ vanishes upon inserting any section of the distribution

$$\mathcal{G} \times_H \mathfrak{p}/\mathfrak{h} \subseteq \mathcal{G} \times_H \mathfrak{g}/\mathfrak{h} \cong T\tilde{M};$$

see, e.g., [11, Proposition 1.5.13]. Conversely, the correspondence space theorem [11, Theorem 1.5.14] states that if $(\mathcal{G} \to \tilde{M}, \omega)$ is a Cartan geometry of type $(G, H)$ and its curvature $\kappa$ satisfies

$$\kappa(\mathcal{G} \times_H \mathfrak{p}/\mathfrak{h}, T\tilde{M}) = 0,$$

then $(\mathcal{G} \to \tilde{M}, \omega)$ is locally isomorphic to the result of the correspondence space construction; see the first part of Definition 2 below.

**Definition 2** Let $(\mathcal{G} \to \tilde{M}, \omega)$ be a Cartan geometry of type $(G, H)$ and $\mathfrak{n} \subseteq \mathfrak{g}$ a $P$-submodule. Denote by $p : \mathcal{G} \to \tilde{M}$ the natural projection. We say that $(\mathcal{G} \to \tilde{M}, \omega)$ is locally isomorphic to the result of the correspondence space construction if there exists an open cover $\tilde{M} = \bigcup_{i\in\mathcal{I}} \tilde{U}_i$ and a collection of Cartan geometries $\{(\mathcal{G}_i \to M_i, \omega_i) : i \in \mathcal{I}\}$ of type $(G, P)$ such that the following holds. Denote by $\{q_i : \mathcal{G}_i \to \mathcal{G}_i/H\}_{i\in\mathcal{I}}$ the natural projections. For every $i \in \mathcal{I}$, there exists an open subset $\tilde{V}_i \subseteq \mathcal{G}_i/H$ and an isomorphism

$$\Phi_i : (p^{-1}(\tilde{U}_i) \to \tilde{U}_i, \omega|_{p^{-1}(\tilde{U}_i)}) \xrightarrow{\cong} (q_i^{-1}(\tilde{V}_i) \to \tilde{V}_i, \omega_i|_{q_i^{-1}(\tilde{V}_i)})$$

of Cartan geometries of type $(G, H)$.

More generally, we say that $(\mathcal{G} \to \tilde{M}, \omega)$ is locally isomorphic to the result of the correspondence space construction modulo $L(\mathfrak{g}/\mathfrak{p}, \mathfrak{n})$ if there exists an open cover $\tilde{M} = \bigcup_{i\in\mathcal{I}} \tilde{U}_i$ and a collection of Cartan geometries $\{(\mathcal{G}_i \to M_i, \omega_i) : i \in \mathcal{I}\}$ of type $(G, P)$ such that the following holds. Denote by $\{q_i : \mathcal{G}_i \to \mathcal{G}_i/H\}_{i\in\mathcal{I}}$ the natural projections. For every $i \in \mathcal{I}$, there exists an open subset $\tilde{V}_i \subseteq \mathcal{G}_i/H$ and an isomorphism

$$\Phi_i : (p^{-1}(\tilde{U}_i) \to \tilde{U}_i) \xrightarrow{\cong} (q_i^{-1}(\tilde{V}_i) \to \tilde{V}_i)$$

of principal $H$-bundles such that $\Phi_i^*\omega_i \equiv \omega$ modulo $L(\mathfrak{g}/\mathfrak{p}, \mathfrak{n})$, i.e.,

$$\begin{cases} \Phi_i^*\omega_i - \omega \in \Omega^1(p^{-1}(\tilde{U}_i), \mathfrak{n}) \quad \text{and} \\ (\omega|_{p^{-1}(\tilde{U}_i)})^{-1}(X) = (\Phi_i^*\omega_i)^{-1}(X) \quad \text{for all } X \in \mathfrak{p}. \end{cases}$$

In [21, Proposition 4.13], the authors have developed a generalization of the classical correspondence space theorem [11, Theorem 1.5.14] recalled just before Definition 2 in their specific setting. The proof of Proposition 3 below expands the general line of their argument. In the special case where $\mathfrak{n} = 0$, it is just the classical correspondence space theorem.

**Proposition 3** *Let $H \subseteq P \subseteq G$ be as at the beginning of this section. Let $\mathfrak{n} \subseteq \mathfrak{g}$ be a $P$-submodule. Given a Cartan geometry $(\mathcal{G} \to \tilde{M}, \omega)$ of type $(G, H)$ with curvature $\kappa$, consider the distribution*

$$\mathcal{V}\tilde{M} = \mathcal{G} \times_H \mathfrak{p}/\mathfrak{h} \subseteq T\tilde{M}.$$

*The following are equivalent.*

*(i) $\kappa(\mathcal{V}\tilde{M}, \mathcal{V}\tilde{M}) = 0$ and $\kappa(\mathcal{V}\tilde{M}, T\tilde{M}) \subseteq \mathcal{G} \times_H \mathfrak{n}$.*
*(ii) $(\mathcal{G} \to \tilde{M}, \omega)$ is locally the result of the correspondence space construction modulo $L(\mathfrak{g}/\mathfrak{p}, \mathfrak{n})$.*

***Proof*** Suppose that $(\mathcal{G} \to M, \omega')$ is a Cartan geometry of type $(G, P)$, $\tilde{M} = \mathcal{G}/H$ and $\omega' \equiv \omega$ modulo $L(\mathfrak{g}/\mathfrak{p}, \mathfrak{n})$; see Definition 2. Then for the curvature $\kappa'$ of $\omega'$, there holds $\kappa'(\mathcal{V}\tilde{M}, T\tilde{M}) = 0$. To prove that (ii) implies (i), it suffices to show that

$$(\kappa' - \kappa)(\omega^{-1}(X), \omega^{-1}(Y)) = 0 \quad \text{and}$$
$$(\kappa' - \kappa)(\omega^{-1}(X), \omega^{-1}(A)) \subseteq \mathfrak{n}$$

for all $X, Y \in \mathfrak{p}$ and $A \in \mathfrak{g}$. Indeed, by assumption, $\omega$ agrees with $\omega'$ on the vertical bundle of $\mathcal{G} \to M$. Since $\mathfrak{n} \subseteq \mathfrak{g}$ is a $P$-submodule, there holds $[\mathfrak{p}, \mathfrak{n}] \subseteq \mathfrak{n}$. The above two relations follow from the formula for curvature in Definition 1 and a brief computation.

Conversely, we prove that (i) implies (ii) in two steps. The first step is precisely the first half of the proof of the classical correspondence space theorem; see, e.g., [11, Theorem 1.5.14]. Indeed, the assumption

$$\kappa(\mathcal{V}\tilde{M}, \mathcal{V}\tilde{M}) = 0$$

means that $\omega^{-1} : \mathfrak{p} \to \mathfrak{X}(\mathcal{G})$ is a Lie algebra homomorphism. Moreover, it implies that $\mathcal{V}\tilde{M} \subseteq T\tilde{M}$ is an involutive distribution and thus around every point of $\tilde{M}$ we may choose a sufficiently small local leaf space

$$\tilde{M} \supseteq \tilde{W} \twoheadrightarrow U$$

such that there exists a global section of $\tilde{W} \twoheadrightarrow U$ and a local section of $p : \mathcal{G} \to \tilde{M}$ defined on $\tilde{W}$. Denote their composition by

$$\sigma : U \to \mathcal{G}.$$

We choose an arbitrary inner product on $\mathfrak{p}$, thus

$$\mathfrak{p} = \mathfrak{h} \oplus \mathfrak{h}^{\perp}.$$

Let $B_{\mathfrak{h}^{\perp}} \subseteq \mathfrak{h}^{\perp}$ and $B_{\mathfrak{h}} \subseteq \mathfrak{h}$ be sufficiently small open balls such that the maps

$$B_{\mathfrak{h}^{\perp}} \times B_{\mathfrak{h}} \to P, \quad (X, Y) \mapsto \exp X \exp Y$$

and

$$B_{\mathfrak{h}^{\perp}} \to P/H, \quad X \mapsto \exp X \bmod H$$

are both diffeomorphisms onto their images. In particular, the group operation $P \times P \to P$ restricts to an injective map

$$\exp B_{\mathfrak{h}^{\perp}} \times H \to P. \tag{3}$$

It follows from Lie's second fundamental theorem that, shrinking $U$ further, the Lie algebra homomorphism $\omega^{-1} : \mathfrak{p} \to \mathfrak{X}(\mathcal{G})$ integrates to a local $P$-action

$$\mathcal{G} \times P \supseteq \sigma(U) \times (\exp B_{\mathfrak{h}^{\perp}} \times \exp B_{\mathfrak{h}}) \xrightarrow{F} p^{-1}(\tilde{W}) \subseteq \mathcal{G}.$$

Observe that the restriction to $(\exp B_{\mathfrak{h}^{\perp}} \times \exp B_{\mathfrak{h}})$ of the left invariant vector field generated by $X \in \mathfrak{p}$ is related to $\omega^{-1}(X)$ by $F$ and that the tangent map of $F$ is bijective at every point. By shrinking $B_{\mathfrak{h}^\perp}$ and $B_{\mathfrak{h}}$, we may assume that $F$ is a diffeomorphism onto its image. Now we restrict $F$ to $\sigma(U) \times \exp B_{\mathfrak{h}^\perp} \times \{e\}$ and extend it to a local $P$-action

$$\tilde{F} : \sigma(U) \times (\exp B_{\mathfrak{h}^\perp})H \to p^{-1}(\tilde{W}) \subseteq \mathcal{G}$$

that is a diffeomorphism onto its image. This is well-defined in view of the injectivity of (3). Now let $\mathcal{G}' = U \times P$. Denote by $q : \mathcal{G}' \to \mathcal{G}'/H$ the natural projection and by $\operatorname{im} \tilde{F}$ the image of $\tilde{F}$. Let

$$\tilde{V} = U \times \frac{(\exp B_{\mathfrak{h}^\perp})H}{H} \subseteq \mathcal{G}'/H \quad \text{and} \quad \tilde{U} = p(\operatorname{im} \tilde{F}) \subseteq \tilde{W} \subseteq \tilde{M}.$$

Note that

$$q^{-1}(\tilde{V}) = U \times (\exp B_{\mathfrak{h}^\perp})H \subseteq \mathcal{G}' \quad \text{and} \quad p^{-1}(\tilde{U}) = \operatorname{im} \tilde{F}.$$

Identifying $\sigma(U)$ with $U$, we obtain an isomorphism

$$\Phi = \tilde{F} : (q^{-1}(\tilde{V}) \to \tilde{V}) \to (p^{-1}(\tilde{U}) \to \tilde{U})$$

of principal $H$-bundles. By construction, $\Phi^*\omega$ is a Cartan connection of type $(G, H)$ that in addition reproduces the generators in $\mathfrak{p}$ of the fundamental vector fields restricted to the subspace $q^{-1}(\tilde{V})$. Hence there exists a unique Cartan connection $\omega' \in \Omega^1(\mathcal{G}', \mathfrak{g})^P$ of type $(G, P)$ such that

$$\omega'|_{U\times H} = \Phi^*\omega|_{U\times H}.$$

In addition, there holds

$$(\omega'|_{q^{-1}(\tilde{V})})^{-1}(X) = (\Phi^*\omega)^{-1}(X) \tag{4}$$

for all $X \in \mathfrak{p}$.

The second step is a minor modification of the second half of the proof of the classical correspondence space theorem. Assume in addition that

$$\kappa(\mathcal{V}\tilde{M}, T\tilde{M}) \subseteq \mathcal{G} \times_H \mathfrak{n}.$$

It remains to show that $\Phi^*\omega - \omega' \in \Omega^1(q^{-1}(\tilde{V}), \mathfrak{n})$. Recall from the first step that $\Phi^*\omega$ coincides with $\omega'$ on $U \times \{e\}$. Moreover, recall that by assumption, every element in $(\exp B_{\mathfrak{h}^\perp})H$ may be expressed as $(\exp X)h$ for a unique choice of $X \in B_{\mathfrak{h}^\perp}$ and $h \in H$. Hence it suffices to verify that

$$(r^{\exp(X)})^*(\Phi^*\omega) - Ad(\exp(-X)) \circ \Phi^*\omega \in \Omega^1(q^{-1}(\tilde{V}), \mathfrak{n})$$

for all $X \in B_{\mathfrak{h}^\perp}$. Indeed, denote by $\mathcal{L}$ the Lie derivative and by the lower dot the directional derivative. Then for all $\xi, \eta \in \mathfrak{X}(\mathcal{G})$, there holds

$$(\mathcal{L}_\xi\omega)(\eta) = \xi.\omega(\eta) - \omega([\xi, \eta])$$

and so

$$\kappa(\xi, \eta) = (\mathcal{L}_\xi \omega)(\eta) - \eta.\omega(\xi) + [\omega(\xi), \omega(\eta)].$$

Substituting $\xi = \omega^{-1}(A)$ for any $A \in \mathfrak{g}$, we obtain

$$\kappa(\omega^{-1}(A), \cdot) = \mathcal{L}_{\omega^{-1}(A)}\omega + ad(A) \circ \omega.$$

This relation together with the assumption $\kappa(\mathcal{V}\tilde{M}, T\tilde{M}) \subseteq \mathcal{G} \times_H \mathfrak{n}$ imply that

$$\mathcal{L}_{\omega^{-1}(X)}\omega + ad(X) \circ \omega \in \Omega^1(\mathcal{G}, \mathfrak{n})$$

for all $X \in B_{\mathfrak{h}^\perp}$. Since $\Phi^*(\omega^{-1}(X))$ is the restriction to $q^{-1}(\tilde{V})$ of the fundamental vector field $\zeta_X \in \mathfrak{X}(\mathcal{G}')$ generated by $X$, the pullback of the above expression integrates to

$$(r^{\exp(tX)})^*(\Phi^*\omega) - Ad(\exp(-tX)) \circ \Phi^*\omega \in \Omega^1(q^{-1}(\tilde{V}), \mathfrak{n})$$

for all $t \in [0, 1]$. Putting $t = 1$, (ii) follows as required. □

Compared to the classical correspondence theorem, this variation requires less information on the curvature of the Cartan geometry of interest and is particularly useful in settings where we are flexible to modify the Cartan connections, as is the case when studying Fefferman-type constructions.

## 3 A Characterization Theorem for Abstract Fefferman-Type Constructions

We follow the original setup of the Fefferman-type functors for Cartan geometries in [5, Subsection 4.5] and establish a characterization result for Cartan geometries arising from this construction. We refer the reader also to [13, Subsection 3.2], which generalizes the original construction to a slightly broader class of functors and establishes a relation between the curvature of the initial Cartan geometry and the curvature of the resulting Cartan geometries; see [13, Subsection 3.3].

In this section, we assume that $G$ and $\tilde{G}$ are Lie groups, $\tilde{P} \subseteq \tilde{G}$ is a closed Lie subgroup and

$$i : G \to \tilde{G}$$

is an injective Lie group homomorphism such that the $G$-orbit of $e\tilde{P}$ in $\tilde{G}/\tilde{P}$ is open. Then for any closed Lie subgroup $P \subseteq G$ containing $H = i^{-1}(\tilde{P})$, $i$ defines a Fefferman-type functor from Cartan geometries of type $(G, P)$ to Cartan geometries of type $(\tilde{G}, \tilde{P})$. We describe the resulting Cartan geometry of the functor applied to a Cartan geometry $(\mathcal{G} \to M, \omega)$ of type $(G, P)$. Let $\tilde{M} = \mathcal{G}/H$ and

$$j : \mathcal{G} \to \mathcal{G} \times_H \tilde{P} = \tilde{\mathcal{G}}$$

be the extension of the principal $H$-bundle $\mathcal{G} \to \tilde{M}$ based on the group homomorphism $i|_H : H \to \tilde{P}$. We obtain $(\tilde{\mathcal{G}} \to \tilde{M}, \tilde{\omega})$ where $\tilde{\omega}$ is the unique Cartan connection of type $(\tilde{G}, \tilde{P})$ on $\tilde{\mathcal{G}} \to \tilde{M}$ such that

$$j^*\tilde{\omega} = i' \circ \omega.$$

Denote the Lie algebras of the Lie groups

$$H \subseteq P \subseteq G \quad \text{and} \quad \tilde{P} \subseteq \tilde{G}$$

by $\mathfrak{h} \subseteq \mathfrak{p} \subseteq \mathfrak{g}$ and $\tilde{\mathfrak{p}} \subseteq \tilde{\mathfrak{g}}$, respectively. Recall that the $G$-orbit of $e\tilde{P}$ in $\tilde{G}/\tilde{P}$ is open if and only if $G$ acts locally transitively on $\tilde{G}/\tilde{P}$, i.e., if $i'(\mathfrak{g}) + \tilde{\mathfrak{p}} = \tilde{\mathfrak{g}}$. Observe that the composition map $\mathfrak{g} \xrightarrow{i'} \tilde{\mathfrak{g}} \twoheadrightarrow \tilde{\mathfrak{g}}/\tilde{\mathfrak{p}}$ is surjective and descends to an isomorphism $\mathfrak{g}/\mathfrak{h} \cong \tilde{\mathfrak{g}}/\tilde{\mathfrak{p}}$ of $H$-modules. We denote by

$$\pi : \tilde{\mathfrak{g}}/\tilde{\mathfrak{p}} \twoheadrightarrow \mathfrak{g}/\mathfrak{p}$$

the composition of the natural projection $\mathfrak{g}/\mathfrak{h} \twoheadrightarrow \mathfrak{g}/\mathfrak{p}$ with the inverse of the above isomorphism. Then the curvature $\tilde{\kappa}$ of $\tilde{\omega}$ is completely determined by the curvature $\kappa$ of $\omega$ via the relation

$$\tilde{\kappa}(j(u)) = i' \circ \kappa(u) \circ \wedge^2 \pi \tag{5}$$

for all $u \in \mathcal{G}$.

Conversely, the subclass of all Cartan geometries of type $(\tilde{G}, \tilde{P})$ that is locally the result of the Fefferman-type functor may be characterized in terms of the existence of a certain parallel tractor. Here, being locally the result of a Fefferman-type functor is defined in a way analogous to Definition 2. Assume now that $(\tilde{\mathcal{G}} \to \tilde{M}, \tilde{\omega})$ is any Cartan geometry of type $(\tilde{G}, \tilde{P})$. Consider the canonical extension $\tilde{\mathcal{G}} \to \tilde{\mathcal{G}} \times_{\tilde{P}} \tilde{G}$ of $\tilde{\mathcal{G}}$ to a principal $\tilde{G}$-bundle. There is a unique principal connection $\tilde{\theta} \in \Omega^1(\tilde{\mathcal{G}} \times_{\tilde{P}} \tilde{G}, \tilde{\mathfrak{g}})$ which pulls back to $\tilde{\omega}$; see, e.g., [11, Theorem 1.5.6]. In particular, given any $\tilde{G}$-representation $\tilde{\mathbb{W}}$, we obtain a tractor bundle

$$\tilde{\mathcal{W}} = (\tilde{\mathcal{G}} \times_{\tilde{P}} \tilde{G}) \times_{\tilde{G}} \tilde{\mathbb{W}} \cong \tilde{\mathcal{G}} \times_{\tilde{P}} \tilde{\mathbb{W}}$$

with tractor connection induced by the principal connection $\tilde{\theta}$.

Assume in addition that there exists $\tilde{w}_0 \in \tilde{\mathbb{W}}$ such that

$$i(G) = \mathrm{Stab}_{\tilde{G}}(\tilde{w}_0) = \{\tilde{g} \in \tilde{G} : \tilde{g}\,\tilde{w}_0 = \tilde{w}_0\}.$$

The $\tilde{P}$-orbit of $\tilde{w}_0$ defines a fiber subbundle

$$\tilde{\mathcal{O}} = \tilde{\mathcal{G}} \times_{\tilde{P}} (\tilde{P}\,\tilde{w}_0) \subseteq \tilde{\mathcal{W}}$$

with standard fiber $\tilde{P}\,\tilde{w}_0 \cong \tilde{P}/i(H)$. The following characterization may be considered as a kin of [9, Theorem 2.6 (ii)] in a simplified setup.

**Theorem 4** *Notation as above. $(\tilde{\mathcal{G}} \to \tilde{M}, \tilde{\omega})$ is locally the result of the Fefferman-type functor applied to a Cartan geometry of type $(G, P)$ if and only if there exists a parallel section*

$$\sigma \in C^\infty(\tilde{\mathcal{G}}, \tilde{P}\,\tilde{w}_0)^{\tilde{P}} = \Gamma(\tilde{\mathcal{O}}) \subseteq \Gamma(\tilde{\mathcal{W}})$$

*of the tractor bundle* $\tilde{\mathcal{W}}$ *such that whenever* $\sigma(u) = \tilde{w}_0$ *where* $u \in \tilde{\mathcal{G}}$, *there holds*

$$\tilde{\kappa}(u)(i'(\mathfrak{p}), \tilde{\mathfrak{g}}) = 0.$$

***Proof*** Let us identify $H \subseteq G$ with their images under $i$ and identify their Lie algebras with their images under $i'$.

Assume that $(\tilde{\mathcal{G}} \to \tilde{M}, \tilde{\omega})$ is the result of the Fefferman-type functor applied to a Cartan geometry $(\mathcal{G} \to M, \omega)$ of type $(G, P)$. Since $H$ fixes $\tilde{w}_0$, the constant map $\sigma : \mathcal{G} \to \{\tilde{w}_0\}$ is $H$-equivariant. It is extended to

$$\tilde{\sigma} \in C^\infty(\tilde{\mathcal{G}}, \tilde{P}\,\tilde{w}_0)^{\tilde{P}} = \Gamma(\tilde{\mathcal{O}})$$

in a unique way. Denote by $p : \mathcal{G} \to \tilde{M}$ the natural projection. For every tangent vector in $T\mathcal{G}$, there holds

$$\nabla^{\tilde{\mathcal{W}}}_{Tp(\,\cdot\,)}\tilde{\sigma} = d\,\sigma(\,\cdot\,) + \omega(\,\cdot\,)\,\sigma. \tag{6}$$

Here, $d\,\sigma$ denotes the directional derivative of the constant map $\sigma : \mathcal{G} \to \{\tilde{w}_0\}$ and $\omega(\,\cdot\,)\,\sigma$ denotes the induced Lie algebra action of $\mathfrak{g}$ on $\tilde{w}_0$. Both of them are zero here. In particular, $\tilde{\sigma}$ is a parallel section of $\tilde{\mathcal{W}}$. Moreover, it follows from the definition of $\tilde{\mathcal{W}}$ that for $u \in \tilde{\mathcal{G}}$, there holds $\tilde{\sigma}(u) = \tilde{w}_0$ if and only if $u \in \mathcal{G}$. In this case, it follows from (5) that $\tilde{\kappa}(u)(i'(\mathfrak{p}), \tilde{\mathfrak{g}}) = 0$.

Conversely, let $\tilde{\sigma} \in C^\infty(\tilde{\mathcal{G}}, \tilde{P}\,\tilde{w}_0)^{\tilde{P}} = \Gamma(\tilde{\mathcal{O}})$ be as stated in the theorem and

$$\mathcal{G} = \{u \in \tilde{\mathcal{G}} : \tilde{\sigma}(u) = \tilde{w}_0\}.$$

Observe that the natural projection $\mathcal{G} \to \tilde{M}$ admits smooth local sections and that the preimage of every point in $\tilde{M}$ is exactly an $H$-orbit of the fiber of $\tilde{\mathcal{G}}$ over the same point. Thus $\mathcal{G}$ is a reduction of $\tilde{\mathcal{G}}$ to a principal $H$-bundle. It remains to show that $\omega = \tilde{\omega}|_{\mathcal{G}}$ has values in $\mathfrak{g}$ so that, in particular, $(\mathcal{G} \to \tilde{M}, \omega)$ is a Cartan geometry of type $(G, H)$. Indeed, $\sigma = \tilde{\sigma}|_{\mathcal{G}}$ is the $\tilde{w}_0$-valued constant map. Since $\mathcal{G} \subseteq \tilde{\mathcal{G}}$, (6) holds where $\nabla^{\tilde{\mathcal{W}}}\tilde{\sigma}$ and $d\,\sigma$ are both zero. In particular, the values of $\omega$ annihilate $\tilde{w}_0$. Since $G = \mathrm{Stab}_{\tilde{G}}(\tilde{W}_0)$, there holds $\mathfrak{g} = \mathrm{Ann}_{\tilde{\mathfrak{g}}}(\tilde{w}_0)$. In particular, the images of $\omega$ lie in $\mathfrak{g}$. This proves our claim. Moreover, it follows from the assumption on $\tilde{\kappa}$ that $(\mathcal{G} \to \tilde{M}, \omega)$ is locally the correspondence space of a Cartan geometry of type $(G, P)$. That is, $(\tilde{\mathcal{G}} \to \tilde{M}, \tilde{\omega})$ is locally the result of the Fefferman-type functor applied to a Cartan geometry of type $(G, P)$. This completes the proof of the theorem. □

One may want to apply the above theorem recursively to obtain the characterization of a certain Fefferman-type construction. For example, denote by $G \subseteq \hat{G} \subseteq \tilde{G}$ the natural inclusions

$$SU(p+1, q+1) \subseteq U(p+1, q+1) \subseteq SO(2p+2, 2q+2).$$

The Fefferman-type construction associating to every CR geometry of type $(p, q)$ a conformal geometry of signature $(2p+1, 2q+1)$ is defined by the canonical group homomorphism $G \to \tilde{G}$. There exists $\mathbb{J} \in \mathfrak{so}(2p+2, 2q+2)$ such that $\mathrm{Stab}_{\tilde{G}}(\mathbb{J}) = \hat{G}$. Moreover, for any $0 \neq \nu \in \wedge_{\mathbb{C}}^{p+q+2} \mathbb{C}^{p+1,q+1}$, there holds $\mathrm{Stab}_{\hat{G}}(\nu) = G$. Thus we arrive at the characterizing condition provided in [8, 2.5 Proposition]. Analogously, the Fefferman-type construction associating to every $(2n+1)$-dimensional torsion-free Lagrangian contact structure of signature $(n+1, n+1)$ is defined by the natural inclusion $SL(n+2, \mathbb{R}) \subseteq SO(n+2, n+2)$. The subgroup of $SO(n+2, n+2)$ fixing a certain point $\mathbb{K} \in \mathfrak{so}(n+2, n+2)$ is $GL(n+2, \mathbb{R})$ and the subgroup of $GL(n+2, \mathbb{R})$ fixing any nonzero element of $\wedge^{n+2} \mathbb{R}$ is $SL(n+2, \mathbb{R})$. This leads to a characterizing condition that is equivalent to [21, Proposition 3.10].

# 4 Some Parabolic Background

Assume now that $G$ is a semisimple Lie group whose Lie algebra $\mathfrak{g}$ comes with a $|k|$-grading, i.e., there is a vector space decomposition $\mathfrak{g} = \bigoplus_{\ell=-k}^{k} \mathfrak{p}_\ell$ such that (i) $[\mathfrak{p}_i, \mathfrak{p}_j] \subseteq \mathfrak{p}_{i+j}$ for all $i, j \in \mathbb{Z}$, (ii) the Lie subalgebra $\mathfrak{p}_- = \bigoplus_{\ell=-k}^{-1} \mathfrak{p}_\ell$ is generated by $\mathfrak{p}_{-1}$ and (iii) $\mathfrak{p}_{\pm k} \neq 0$. We follow the convention that $\mathfrak{p}_i = 0$ for $|i| > k$ and $\mathfrak{p}^i = \bigoplus_{\ell=i}^{k} \mathfrak{p}_\ell$.

Let $\mathfrak{p} = \mathfrak{p}^0$ and $P \subseteq G$ be a Lie subgroup with Lie algebra $\mathfrak{p}$. $P$ is said to be a parabolic subgroup of $G$. Note that the maximal parabolic subgroup of $G$ is

$$\{g \in G : Ad(g)(\mathfrak{p}^i) \subseteq \mathfrak{p}^i, -k \leq i \leq k\}.$$

The Levi subgroup of $P$ is given by

$$P_0 = \{g \in P : Ad(g)(\mathfrak{p}_i) \subseteq \mathfrak{p}_i, -k \leq i \leq k\}.$$

It is the maximal subgroup of $P$ with Lie algebra $\mathfrak{p}_0$. Moreover, let $\mathfrak{p}_+ = \mathfrak{p}^1$ and $P_+ = \exp \mathfrak{p}_+$. Then $P_+$ is a normal subgroup of $P$ and $P/P_+ \cong P_0$.

The Killing form of $\mathfrak{g}$ induces a $G$-invariant identification $\mathfrak{g} \cong \mathfrak{g}^*$, which further induces a $P$-module identification

$$(\mathfrak{g}/\mathfrak{p})^* \cong \mathfrak{p}_+.$$

With this identification, the Kostant codifferential

$$\partial^* : \wedge^{n+1} (\mathfrak{g}/\mathfrak{p})^* \otimes \mathfrak{g} \to \wedge^{n} (\mathfrak{g}/\mathfrak{p})^* \otimes \mathfrak{g} \tag{7}$$

is determined by

$$\begin{aligned}\partial^*(Z_0 \wedge \cdots \wedge Z_n \otimes A) = &\sum_{i=0}^{n} (-1)^{i+1} Z_0 \wedge \cdots \wedge \hat{Z}_i \wedge \cdots \wedge Z_n \otimes [Z_i, A] \\ &+ \sum_{i<j} (-1)^{i+j} [Z_i, Z_j] \wedge Z_0 \wedge \cdots \wedge \hat{Z}_i \wedge \cdots \wedge \hat{Z}_j \wedge \cdots \wedge Z_n \otimes A\end{aligned}$$

for $Z_0, \dots, Z_n \in \mathfrak{p}_+$, $A \in \mathfrak{g}$; see, e.g., [11, Subsection 3.1.11]. Here, the hat indicates that the term is omitted from a series. Note that $\partial^*$ is $P$-equivariant and $\partial^* \circ \partial^* = 0$.

In particular, let $\{X^i\} \subseteq \mathfrak{g}/\mathfrak{p}$ be a basis with dual basis $\{Z_i\} \subseteq \mathfrak{p}_+$. For $\phi \in \wedge^2 (\mathfrak{g}/\mathfrak{p})^* \otimes \mathfrak{g}$ and $X \in \mathfrak{g}/\mathfrak{p}$, we have

$$(\partial^*\phi)(X) = 2\sum_i [Z_i, \phi(X, X^i)] - \sum_i \phi([Z_i, \tilde{X}], X^i) \tag{8}$$

where $\tilde{X} \in \mathfrak{g}$ is any lift of $X \in \mathfrak{g}/\mathfrak{p}$; see, e.g., [11, Lemma 3.1.11]. Note that if $\mathfrak{g}$ is $|1|$-graded, the second summand on the right-hand side vanishes as $[\mathfrak{g}, \mathfrak{p}_+] \subseteq \mathfrak{p}$.

Since $P \subseteq G$ is a parabolic subgroup, a Cartan geometry $(\mathcal{G} \to M, \omega)$ of type $(G, P)$ is called a parabolic geometry of type $(G, P)$. It is said to be regular if its curvature $\kappa : \mathcal{G} \to L(\wedge^2 \mathfrak{g}/\mathfrak{p}, \mathfrak{g})$ has homogeneity $\geq 1$, i.e., if $\kappa(u)(\mathfrak{p}^i, \mathfrak{p}^j) \subseteq \mathfrak{p}^{i+j+1}$ for all $u \in \mathcal{G}$. Note that this property holds automatically if $\mathfrak{g}$ is $|1|$-graded. On the other hand, $(\mathcal{G}, \omega)$ is said to be normal if $\partial^*\kappa = 0$, i.e.

$$\kappa : \mathcal{G} \to \ker\, \partial^* \subseteq L(\wedge^2 \mathfrak{g}/\mathfrak{p}, \mathfrak{g}).$$

In this case, $\kappa$ projects to

$$\kappa_H : \mathcal{G} \to \ker\, \partial^*/\operatorname{im} \partial^*.$$

This projection is called the harmonic curvature of $(\mathcal{G} \to M, \omega)$.

The Lie algebra homology differential

$$\partial : L(\wedge^n \mathfrak{p}_-, \mathfrak{g}) \to L(\wedge^{n+1} \mathfrak{p}_-, \mathfrak{g}) \tag{9}$$

given by

$$\begin{aligned}(\partial\phi)(X^0, \dots, X^n) &= \sum_{i=0}^{n} (-1)^i\, [X^i, \phi(X^0, \dots, \hat{X}^i, \dots, X^n)] \\ &\quad + \sum_{i<j} \phi([X^i, X^j], X^0, \dots, \hat{X}^i, \dots, \hat{X}^j, \dots, X^n)\end{aligned}$$

is $P_0$-equivariant; see, e.g., [11, Subsection 3.1.10]. Here, $X^0, \dots, X^n \in \mathfrak{p}_-$ and $\phi \in L(\wedge^n \mathfrak{p}_-, \mathfrak{g})$. The hat indicates that the object is omitted. Then $\Box = \partial \circ \partial^* + \partial^* \circ \partial$ is the Kostant Laplacian; see, e.g., [11, Subsection 3.1.11]. Recall from [11, Proposition 3.1.11] the $P_0$-invariant Hodge decomposition

$$L(\wedge^2 \mathfrak{g}/\mathfrak{p}, \mathfrak{g}) = \operatorname{im} \partial^* \oplus \ker \Box \oplus \operatorname{im} \partial \tag{10}$$

where in addition, $\operatorname{im} \partial^* \oplus \ker \Box = \ker\, \partial^*$ and $\ker \Box \oplus \operatorname{im} \partial = \ker\, \partial$. In particular, one obtains an identification

$$\kappa_H(u) \in \ker\, \partial^*/\operatorname{im} \partial^* \cong \ker \Box \subseteq L(\wedge^2 \mathfrak{g}/\mathfrak{p}, \mathfrak{g})$$

for every $u \in \mathcal{G}$. Using this, one can rephrase the improved Bianchi identity as follows.

**Lemma 5** *(See [12, Proposition 4.16 (3)]) Let $(\mathcal{G}, \omega)$ be a normal regular parabolic geometry of type $(G, P)$ with curvature*

$$\kappa : \mathcal{G} \to \ker \partial^* \subseteq L(\wedge^2 \mathfrak{g}/\mathfrak{p}, \mathfrak{g})$$

*and harmonic curvature $\kappa_H$. Let $\mathbb{F}$ be a $P$-submodule of $L(\wedge^2 \mathfrak{g}/\mathfrak{p}, \mathfrak{g})$ such that, viewing $\kappa_H(u) \in \ker \Box$, we have*

$$\kappa_H(u) \in \mathbb{F} \cap \ker \Box$$

*for all $u \in \mathcal{G}$ and such that $\mathbb{F}$ is stable under $\mathbb{F}$-insertions, i.e., for every $\varphi, \psi \in \mathbb{F}$, applying the Kostant codifferential $\partial^*$ to the map*

$$(X, Y, Z) \mapsto \psi(\varphi(X, Y), Z) + \psi(\varphi(Y, Z), X) + \psi(\varphi(Z, X), Y),$$

*where $X, Y, Z \in \mathfrak{g}/\mathfrak{p}$, one obtains an element in $\mathbb{F}$.*

*Then $\kappa(u) \in \mathbb{F}$ for all $u \in \mathcal{G}$.*

**Convention 6** *From now on through the rest of the paper, we fix the notation for the following Lie groups and Lie algebras. Let*

$$H \subseteq Q \subseteq P \subseteq G$$

*be the groups*

$$\left\{\left(\begin{array}{c|c|c} * & * & * \\ \hline 0 & 1 & * \\ \hline 0 & 0 & * \end{array}\right)\right\} \subseteq \left\{\left(\begin{array}{c|c|c} * & * & * \\ \hline 0 & \mathbb{R}_+ & * \\ \hline 0 & 0 & * \end{array}\right)\right\} \subseteq \left\{\left(\begin{array}{c|c|c} * & * & * \\ * & * & * \\ \hline 0 & 0 & * \end{array}\right)\right\} \subseteq SL(n+2, \mathbb{R}),$$

*where the block size is $(1, 1, n) \times (1, 1, n)$. Their Lie algebras are denoted by*

$$\mathfrak{h} \subseteq \mathfrak{q} \subseteq \mathfrak{p} \subseteq \mathfrak{g},$$

*respectively. Note that $Q$ and $P$ are parabolic subgroups of the semisimple Lie group $G$ arising from the grading decompositions*

$$\mathfrak{g} = \left(\begin{array}{c|c|c} \mathfrak{q}_0 & \mathfrak{q}_1^E & \mathfrak{q}_2 \\ \hline \mathfrak{q}_{-1}^E & \mathfrak{q}_0 & \mathfrak{q}_1^V \\ \hline \mathfrak{q}_{-2} & \mathfrak{q}_{-1}^V & \mathfrak{q}_0 \end{array}\right) \quad \text{and} \quad \mathfrak{g} = \left(\begin{array}{c|c} \mathfrak{p}_0 & \mathfrak{p}_1 \\ \hline \mathfrak{p}_{-1} & \mathfrak{p}_0 \end{array}\right),$$

*respectively. We also fix the notation for*

$$\tilde{P} \subseteq \tilde{G} = SL(n+3, \mathbb{R})$$

*where* $\tilde{P} \subseteq \tilde{G}$ *is given by* $P \subseteq G$ *with* $n$ *increased by* 1. *Their Lie algebras are* $\tilde{\mathfrak{p}} \subseteq \tilde{\mathfrak{g}}$, *respectively. The grading*

$$\tilde{\mathfrak{g}} = \tilde{\mathfrak{p}}_{-1} \oplus \tilde{\mathfrak{p}}_0 \oplus \tilde{\mathfrak{p}}_1.$$

*is precisely the grading of* $(\mathfrak{g}, \mathfrak{p})$ *with* $n$ *increased by* 1.

## 4.1 Almost Grassmann Structures

The following categorical correspondence follows from, e.g., [11, Subsection 4.1.3 and Theorem 3.1.14].

**Lemma 7** *Notation as in Convention 6. Every parabolic geometry* $(\mathcal{G} \to M, \omega)$ *of type* $(G, P)$ *has an underlying almost Grassmann structure*

$$(M, E = \mathcal{G} \times_P \mathbb{R}^2, F = \mathcal{G} \times_P \mathbb{R}^{n+2}/\mathbb{R}^2)$$

*of type* $(2, n)$. *Conversely, for every almost Grassmann structure* $(M, E, F)$ *of type* $(2, n)$, *there is, up to isomorphisms, a unique normal parabolic geometry of type* $(G, P)$ *whose underlying geometry is* $(M, E, F)$.

Let $(\mathcal{G} \to M, \omega)$ be the normal parabolic geometry of type $(G, P)$ canonically associated to a given almost Grassmann structure $(M, E, F)$. Recall from, e.g., [11, Step (E) in Subsection 4.1.3] that for $n > 2$, $\kappa_H$ decomposes into

$$\begin{aligned} &\tau \in \Gamma\left((\mathrm{Sym}^2 E \otimes E^*)_o \otimes (\wedge^2 F^* \otimes F)_o\right) \subseteq \Omega^2(M, TM) \quad \text{and} \\ &\rho \in \Gamma(\wedge^2 E \otimes \mathrm{Sym}^2 F^* \otimes \mathfrak{sl}(F)) \subseteq \Omega^2(M, \mathfrak{sl}(F)). \end{aligned} \tag{11}$$

Here, the natural decomposition

$$\wedge^2 T^*M = \wedge^2 (E \otimes F^*) = (\mathrm{Sym}^2 E \otimes \wedge^2 F^*) \oplus (\wedge^2 E \otimes \mathrm{Sym}^2 F^*)$$

is used and

$$(\mathrm{Sym}^2 E \otimes E^*)_o \subseteq \mathrm{Sym}^2 E \otimes E^* \quad \text{and} \quad (\wedge^2 F^* \otimes F)_o \subseteq \wedge^2 F^* \otimes F$$

denote the trace-free component. Stretching notation, we call $\tau$ the harmonic torsion of $(M, E, F)$ and $\rho$ the harmonic curvature of $(M, E, F)$. On the other hand, recall from, e.g., [11, Theorem 3.1.12] that the projection of the curvature $\kappa$ of $\omega$ to

$$T \in \Omega^2_{hor}(\mathcal{G}, \mathfrak{g}/\mathfrak{p})^P \cong \Omega^2(M, TM) \tag{12}$$

coincides with $\tau$ and is called the torsion of $(M, E, F)$. If $\tau = 0$, we say that $(M, E, F)$ is torsion-free.

Although Weyl structures and their induced objects are defined for every parabolic geometry, see, e.g., [11, Subsections 5.1–5.2], we only recall them in the setting of almost Grassmann structure. The Lie algebra grading induces the $P_0$-equivariant decompositions

$$\omega = \omega_{\mathfrak{p}_{-1}} + \omega_{\mathfrak{p}_0} + \omega_{\mathfrak{p}_1}$$
$$\kappa = \kappa_{\mathfrak{p}_{-1}} + \kappa_{\mathfrak{p}_0} + \kappa_{\mathfrak{p}_1}$$

of the Cartan connection $\omega$ as well as its curvature $\kappa$.

Observe that $\mathcal{G}_0 = \mathcal{G}/P_+$ is a principal bundle with structure group $P_0 = P/P_+$. Recall from, e.g., [11, Proposition 3.1.15 (2)] that $\omega_{\mathfrak{p}_{-1}}$ descends to a one-form $\theta \in \Omega^1(\mathcal{G}_0, \mathfrak{p}_{-1})^{P_0}$. In fact, $(\mathcal{G}_0, \theta)$ is the $G$-structure of $(M, E, F)$. In particular, the natural projection $\mathcal{G} \twoheadrightarrow \mathcal{G}_0$ induces identifications

$$\begin{aligned} E &\cong \mathcal{G} \times_P \mathbb{R}^2 \cong \mathcal{G}_0 \times_{P_0} \mathbb{R}^2 \\ F &\cong \mathcal{G} \times_P \mathbb{R}^{n+2}/\mathbb{R}^2 \cong \mathcal{G}_0 \times_{P_0} \mathbb{R}^n \\ TM &\cong \mathcal{G} \times_P \mathfrak{g}/\mathfrak{p} \cong \mathcal{G}_0 \times_{P_0} \mathfrak{p}_{-1} \\ T^*M &\cong \mathcal{G} \times_P \mathfrak{p}_1 \cong \mathcal{G}_0 \times_{P_0} \mathfrak{p}_1 \end{aligned}$$

so that

$$\mathcal{G}_0 \times_{P_0} \mathfrak{p}_0 \cong \mathfrak{s}(L(E, E) \oplus L(F, F))$$

where

$$\mathfrak{s}(L(E, E) \oplus L(F, F)) = \Big(L(E, E) \oplus L(F, F)\Big) \cap \mathfrak{sl}(E \oplus F). \tag{13}$$

Recall from, e.g., [11, Subsection 5.1.1] that there exist $P_0$-equivariant sections of the natural projection $\mathcal{G} \twoheadrightarrow \mathcal{G}_0$. We call such a section $\sigma : \mathcal{G}_0 \to \mathcal{G}$ a Weyl structure of $(M, E, F)$. Observe that $\sigma^*\omega_{\mathfrak{p}_{-1}}$ coincides with the one-form $\theta$ defining the $G$-structure and that $\sigma^*\kappa_{\mathfrak{p}_{-1}}$ coincides with the torsion (12). Hence they are independent of the choice of Weyl structures. We obtain the following objects associated to a Weyl structure $\sigma$; see, e.g., [11, Subsections 5.1.2 and 5.2.3].

$$\begin{aligned} \sigma^*\omega_{\mathfrak{p}_0} &\in \Omega^1(\mathcal{G}_0, \mathfrak{p}_0)^{P_0} \text{ (the Weyl connection)} \\ \mathrm{P} = \sigma^*\omega_{\mathfrak{p}_1} &\in \Omega^1_{hor}(\mathcal{G}_0, \mathfrak{p}_1)^{P_0} \cong \Omega^1(M, T^*M) \text{ (the Rho tensor)} \\ (W, W') = \sigma^*\kappa_{\mathfrak{p}_0} &\in \Omega^2_{hor}(\mathcal{G}_0, \mathfrak{p}_0)^{P_0} \\ &\cong \Omega^2(M, \mathfrak{s}(L(E, E) \oplus L(F, F))) \text{ (the Weyl tensor)} \\ Y = \sigma^*\kappa_{\mathfrak{p}_1} &\in \Omega^2_{hor}(\mathcal{G}_0, \mathfrak{p}_1)^{P_0} \cong \Omega^2(M, T^*M) \text{ (the Cotton-York tensor).} \end{aligned} \tag{14}$$

The Weyl connection $\sigma^*\omega_{\mathfrak{p}_0}$ is a principal connection on $\mathcal{G}_0$ and the Weyl tensor $(W, W')$ is decomposed into

$$W \in \Omega^2(M, L(E, E)) \quad \text{and} \quad W' \in \Omega^2(M, L(F, F)).$$

Following the convention of Penrose abstract index notation, see, e.g., [11, Subsection 4.1.3], we use unprimed upper indices for $E$ and lower indices for $E^*$ and primed upper indices for $F$ and lower indices for $F^*$. Recall that a contraction is denoted by writing the same index variable in the two slots that are being contracted. Thus

$$W = W^{A\ B\ C}_{A'B'D}, \quad W' = W'^{A\ B\ C'}_{A'B'D'} \quad \text{and} \quad \tau = \tau^{A\ B\ C'}_{A'B'D}$$

where $\tau$ is the harmonic torsion of $(M, E, F)$.

The following lemma has been stated in [22, Equation 25]. For an explicit proof, we refer the reader to [18, Proposition 4].

**Lemma 8** *We specify the following contractions.*

$$tr(W)^{A\ B}_{A'B'} = W^{A\ I\ B}_{A'B'I} \quad tr(W')^{A\ B}_{A'B'} = W'^{A\ B\ I'}_{A'I'B'} \quad \textit{and} \quad tr(\iota_\tau \tau)^{A\ B}_{A'B'} = \tau^{I\ A\ J'}_{I'A'J}\,\tau^{J\ B\ I'}_{J'B'I}\,.$$

*Then*

$$tr(W) = tr(W') \in \Gamma(\mathrm{Sym}^2 T^*M)$$

*with*

$$tr(\iota_\tau\tau)^{(A\ B)}_{(A'B')} = n\,tr(W)^{(A\ B)}_{(A'B')} \qquad tr(\iota_\tau\tau)^{[A\ B]}_{[A'B']} = (n+4)\,tr(W)^{[A\ B]}_{[A'B']}$$
$$tr(\iota_\tau\tau)^{(A\ B)}_{[A'B']} = tr(W)^{(A\ B)}_{[A'B']} = 0 \qquad tr(\iota_\tau\tau)^{[A\ B]}_{(A'B')} = tr(W)^{[A\ B]}_{(A'B')} = 0.$$

*Here, $[\cdot, \cdot]$ and $(\cdot, \cdot)$ are the alternation and the symmetrization of indices, respectively.*

Recall from, e.g., [11, Theorem 5.2.3 (1)] that Weyl connections of $(M, E, F)$ are exactly the principal connections on $\mathcal{G}_0$ which induce a connection on $TM \cong E^* \otimes F$ with $\partial^*$-closed torsion. Note that the principal connections on $\mathcal{G}_0$ are equivalent to all pairs $(\nabla^E, \nabla^F)$ of connections on $(E, F)$ compatible with the defining isomorphism $\wedge^2 E^* \cong \wedge^n F$. Moreover, since $(\mathcal{G} \to M, \omega)$ is $|1|$-graded, there is a one-to-one correspondence between Weyl structures and Weyl connections of $(\mathcal{G} \to M, \omega)$; see, e.g., [11, Proposition 5.1.1 and Proposition 5.1.6]. Thus we derive the following one-to-one correspondence between Weyl structures and the admissible pairs $(\nabla^E, \nabla^F)$ mentioned above.

**Lemma 9** *Let $(M, E, F)$ be an almost Grassmann structure of type $(2, n)$. Given a Weyl structure, its Weyl connection $(\nabla^E, \nabla^F)$ on $(E, F)$ satisfies the following properties.*

- *$(\nabla^E, \nabla^F)$ is compatible with the defining isomorphism $\wedge^2 E^* \cong \wedge^n F$ and*
- *$(\nabla^E, \nabla^F)$ induces a connection on $TM \cong E^* \otimes F$ whose torsion is exactly the harmonic torsion $\tau$ of $(M, E, F)$; see (11).*

*Conversely, given any pair of linear connections $(\nabla^E, \nabla^F)$ on $(E, F)$ satisfying these two properties, there is a unique Weyl structure of $(M, E, F)$ whose Weyl connection on $(E, F)$ is $(\nabla^E, \nabla^F)$.*

We may use a Weyl structure and its induced objects to describe tractor connections; see Sect. 3. Let $\sigma : \mathcal{G}_0 \to \mathcal{G}$ be a Weyl structure of $(\mathcal{G} \to M, \omega)$. First, it decomposes a tractor bundle into simpler bundles. For example, the adjoint representation of $G$ restricted to $P_0$ admits decomposition $\mathfrak{g} = \mathfrak{p}_{-1} \oplus \mathfrak{p}_0 \oplus \mathfrak{p}_1$. Applying the functor $\mathcal{G}_0 \times_{P_0} \cdot$ to this decomposition results in

$$\mathcal{A}M \overset{\sigma}{\cong} \mathcal{G}_0 \times_{P_0} \mathfrak{g} = TM \oplus \mathfrak{s}(L(E, E) \oplus L(F, F)) \oplus T^*M; \tag{15}$$

see (13). Similarly, for the standard tractor bundle induced from the standard representation of $G$ on $\mathbb{R}^{n+2}$, the $P_0$-module decomposition $\mathbb{R}^{n+2} = \mathbb{R}^2 \oplus \mathbb{R}^n$ induces

$$\mathcal{T} \overset{\sigma}{\cong} \mathcal{G}_0 \times_{P_0} \mathbb{R}^{n+2} = E \oplus F. \tag{16}$$

Second, consider the the horizontal lift $\mathfrak{X}(M) \ni \xi \mapsto \xi^h \in \mathfrak{X}(\mathcal{G}_0)$ of the Weyl connection $\nabla$ as well as the Rho tensor P. Consider the tractor bundle $\mathcal{W} = \mathcal{G} \times_P W$ induced by a representation $\rho : G \to GL(W)$ with derivative $\rho' : \mathfrak{g} \to L(W, W)$. Then the tractor connection $\nabla^{\mathcal{W}}$ on $\mathcal{W}$ can be expressed as

$$\begin{aligned} \nabla^{\mathcal{W}}_{\xi} s &= \xi^h . f_s + \rho'(\sigma^*\omega_{\mathfrak{p}_{-1}}(\xi^h))(f_s) + \rho'(\sigma^*\omega_{\mathfrak{p}_1}(\xi^h))(f_s) \\ &= \nabla_{\xi} s + \xi \bullet s + \mathrm{P}(\xi) \bullet s; \end{aligned} \tag{17}$$

see, e.g., [11, Proposition 5.1.10 (2)]. Here, $s \in \Gamma(\mathcal{W})$ corresponds to $f_s \in C^\infty(\mathcal{G}_0, W)^{P_0}$, $\xi \in \mathfrak{X}(M)$ and $\bullet : TM \times \mathcal{W} \to \mathcal{W}$ and $\bullet : T^*M \times \mathcal{W} \to \mathcal{W}$ are the $\mathcal{G}_0$-associated maps of $\rho'|_{\mathfrak{p}_{-1}}$ and $\rho'|_{\mathfrak{p}_1}$, respectively.

## 4.2 Path Geometries

Let $\mathcal{D}$ be a distribution on a manifold $M$. The Levi bracket

$$\mathcal{L} : \mathcal{D} \times \mathcal{D} \to TM/\mathcal{D}$$

is the skew-symmetric $C^\infty(M)$-bilinear map descended from the Lie bracket of sections of $\mathcal{D}$. Assume now that $M^{2n+1}$ is $(2n+1)$-dimensional. Recall that a generalized path geometry on $M$ is a distribution $E \oplus V \subseteq TM$ on $M$ where $E$ and $V$ are of rank 1 and rank $n$, respectively, such that the Levi bracket

$$\mathcal{L} : \wedge^2(E \oplus V) \to TM/(E \oplus V)$$

has kernel $\wedge^2 V$ and hence induces an isomorphism

$$E \otimes V \underset{\cong}{\overset{\mathcal{L}}{\to}} TM/(E \oplus V);$$

see, e.g., [11, Definition 4.4.3].

Notice that an orientation on $E \oplus V$ may be encoded as an orientation on $E^* \otimes \wedge^n V$ via the canonical orientation on $\otimes^2 E$.

**Lemma 10** *Notation as in Convention 6. Let $(\mathcal{G} \to M, \omega)$ be a regular parabolic geometry of type $(G, Q)$. Then $M$ has dimension $2n+1$ and*

$$(M, E = \mathcal{G} \times_Q (\mathfrak{q}_{-1}^E \oplus \mathfrak{q})/\mathfrak{q}, V = \mathcal{G} \times_Q (\mathfrak{q}_{-1}^V \oplus \mathfrak{q})/\mathfrak{q})$$

*is a generalized path geometry with oriented rank-$(n+1)$-distribution $E \oplus V \subseteq TM$ or, equivalently, with oriented line bundle $E^* \otimes \wedge^n V$.*

*Conversely, for any generalized path geometry on a $(2n+1)$-dimensional manifold with oriented rank-$(n+1)$-distribution or, equivalently, with oriented line bundle $E^* \otimes \wedge^n V$, there is a unique normal regular parabolic geometry $(\mathcal{G} \to M, \omega)$ of type $(G, Q)$ with such an underlying geometry.*

***Proof*** Denote by $Aut_{gr}(\mathfrak{q}_-)$ the group of Lie algebra automorphisms on $\mathfrak{q}_- = \mathfrak{q}_{-2} \oplus \mathfrak{q}_{-1}$ preserving the grading decomposition. Using this identification, one observes that the adjoint representation restricted to the Levi subgroup $Q_0$ of $Q$ induces an identification

$$Q_0 \xrightarrow[\cong]{Ad} \{\varphi \in Aut_{gr}(\mathfrak{q}_-) : \varphi(\mathfrak{q}_{-1}^E) \subseteq \mathfrak{q}_{-1}^E, \varphi(\mathfrak{q}_{-1}^V) \subseteq \mathfrak{q}_{-1}^V \text{ and } \varphi|_{\mathfrak{q}_{-1}} \in GL_+(\mathfrak{q}_{-1})\}.$$

The fact that $\varphi|_{\mathfrak{q}_{-1}}$ is orientation-preserving is the only difference to the case dealt in [11, Proposition 4.4.3], so the proof can be completed as there. □

Let $(\mathcal{G} \to M, \omega)$ be the normal regular parabolic geometry of type $(G, Q)$ canonically associated to a given path geometry $(M, E, V)$ with oriented $E \oplus V$. Denote by $\kappa_H$ its harmonic curvature. Since $V$ is involutive for a path geometry $(M, E, V)$, recall from, e.g., [11, Subsection 4.4.3] that for $n \geq 2$, $\kappa_H$ decomposes to

$$\begin{aligned} &\tau \in \Gamma(E^* \otimes (TM/(E \oplus V))^* \otimes V) \quad \text{and} \\ &\rho \in \Gamma(V^* \otimes (TM/(E \oplus V))^* \otimes L(V, V)). \end{aligned} \tag{18}$$

Stretching notation, we call $\tau$ the harmonic torsion of $(M, E, V)$ and $\rho$ the harmonic curvature of $(M, E, V)$. Analogous to the setting of almost Grassmann structure, $\tau$ vanishes if and only if the torsion $T$ of $(\mathcal{G} \to M, \omega)$ vanishes. In this case we say that the path geometry $(M, E, V)$ is torsion-free. Here, $T$ is defined by the projection of $\kappa \in \Omega^2_{hor}(\mathcal{G}, \mathfrak{g})^Q$ to

$$T \in \Omega^2_{hor}(\mathcal{G}, \mathfrak{g}/\mathfrak{q})^Q \cong \Omega^2(M, TM).$$

## 5 The Fefferman-Type Construction on Path Geometries

We follow Convention 6 for the notation of Lie groups, Lie algebras and gradings. Denote by $e_1, \ldots, e_{n+3}$ the standard basis of $\mathbb{R}^{n+3}$. The decomposition of $\mathbb{R}^{n+3}$ into the line $\mathbb{R}e_3$ and the hyperplane

$$\mathrm{span}\{e_1, e_2 + e_3, e_4, \dots, e_{n+3}\}$$

induces a representation of the Lie group $G$ on $\mathbb{R}^{n+3}$. We write this representation as the injective group homomorphism

$$i : G \to \tilde{G}, \tag{19}$$

$$\left(\begin{array}{c|c|c} x_{11} & x_{12} & x_{13} \\ \hline x_{21} & x_{22} & x_{23} \\ \hline & & \\ x_{31} & x_{32} & x_{33} \\ & & \end{array}\right) \mapsto \left(\begin{array}{c|c|c|c} x_{11} & x_{12} & 0 & x_{13} \\ \hline x_{21} & x_{22} & 0 & x_{23} \\ \hline x_{21} & x_{22} - 1 & 1 & x_{23} \\ \hline & & & \\ x_{31} & x_{32} & 0 & x_{33} \\ & & & \end{array}\right)$$

where the block sizes are $(1, 1, n) \times (1, 1, n)$ and $(1, 1, 1, n) \times (1, 1, 1, n)$, respectively. Observe that $G$ acts locally transitively on $\tilde{G}/\tilde{P}$, that $H = i^{-1}(\tilde{P})$ and that $H \subseteq Q$. Recall that every path geometry

$$(M, E, V)$$

where $M$ is $(2n + 1)$-dimensional and $E \oplus V$ is oriented is canonically associated to a normal regular parabolic geometry

$$(\mathcal{G} \to M, \omega)$$

of type $(G, Q)$. The Fefferman-type functor defined by $i : G \to \tilde{G}$ associates to $(\mathcal{G} \to M, \omega)$ a parabolic geometry

$$(\tilde{\mathcal{G}} \to \tilde{M}, \tilde{\omega})$$

of type $(\tilde{G}, \tilde{P})$ with underlying almost Grassmann structure

$$(\tilde{M}, \tilde{E}, \tilde{F})$$

of type $(2, n + 1)$. The construction

$$(M, E, V) \mapsto (\tilde{M}, \tilde{E}, \tilde{F}) \tag{20}$$

is clearly functorial.

We denote by

$$p : \tilde{M} \to M$$

the natural projection and by $\kappa$ and $\tilde{\kappa}$ the curvatures of $\omega$ and $\tilde{\omega}$, respectively. Denote by

$$j : \mathcal{G} \to \tilde{\mathcal{G}}$$

the natural map defining $\tilde{\mathcal{G}}$ and by

$$\pi : \tilde{\mathfrak{g}}/\tilde{\mathfrak{p}} \twoheadrightarrow \mathfrak{g}/\mathfrak{q}$$

the projection canonically induced by the Lie algebra map $i'$ of (19); see Sect. 3.

For the statement of the following theorem, recall Lemma 10.

**Theorem 11** *Notation as above. Then*

$$\tilde{M} = GL_+(E^* \otimes \wedge^n V)$$

*as a principal* $\mathbb{R}_+$*-bundle over* $M$ *and*

$$\tilde{E}^* = (Tp)^{-1}(E) \quad \text{and} \quad \tilde{F} = (Tp)^{-1}(V) \tag{21}$$

*as vector bundles over* $\tilde{M}$.

*Moreover, the construction* (20) *is independent of the choice of orientation on* $E^* \otimes \wedge^n V$. *In particular, it leads to a construction associating to every path geometry* $(M', E', V')$ *without auxiliary orientation an almost Grassmann structure* $(\tilde{M}', \tilde{E}', \tilde{F}')$ *where* $\tilde{M}' = GL(E'^* \otimes \wedge^n V')/\mathbb{Z}_2$ *and* $\tilde{E}'^*, \tilde{F}' \subseteq T\tilde{M}'$ *are preimages of* $E', V' \subseteq TM'$, *respectively. Here,* $\mathbb{Z}_2$ *denotes the group generated by the negation map* $v \mapsto -v$ *on the line bundle* $E'^* \otimes \wedge^n V'$.

***Proof*** Observe that the adjoint representation induces a surjective group homomorphism

$$Q \twoheadrightarrow GL_+\left(\left((\mathfrak{q}_{-1}^E \oplus \mathfrak{q})/\mathfrak{q}\right)^* \otimes \wedge^n \left((\mathfrak{q}_{-1}^V \oplus \mathfrak{q})/\mathfrak{q}\right)\right) = \mathbb{R}_+$$

with kernel $H$. Hence

$$\tilde{M} = \mathcal{G}/H = \mathcal{G} \times_H Q/H = GL_+(E^* \otimes \wedge^n V)$$

as a principal $\mathbb{R}_+$-bundle over $M$.

Next, we construct the identifications in (21). Notice that

$$\begin{cases} (Tp)^{-1}(E) = \mathcal{G} \times_H (\mathfrak{q}_{-1}^E \oplus \mathfrak{q})/\mathfrak{h} \quad \text{and} \\ (Tp)^{-1}(V) = \mathcal{G} \times_H (\mathfrak{q}_{-1}^V \oplus \mathfrak{q})/\mathfrak{h}. \end{cases} \tag{22}$$

On the other hand, consider the filtration $\mathbb{R}^{n+2} \supseteq \mathbb{R}^2 \supseteq \mathbb{R}$ induced by the standard representation of $H$ and the filtration $\mathbb{R}^{n+3} \supseteq \mathbb{R}^2$ induced by the standard representation of $\tilde{P}$. Observe that the principal bundle map $j : \mathcal{G} \to \tilde{\mathcal{G}}$ in the Fefferman-type functor canonically induces

$$\begin{cases} \tilde{E}^* = \tilde{\mathcal{G}} \times_{\tilde{P}} \mathbb{R}^{2*} \qquad = \mathcal{G} \times_H \mathbb{R}^{2*} \quad \text{and} \\ \tilde{F} = \tilde{\mathcal{G}} \times_{\tilde{P}} \mathbb{R}^{n+3}/\mathbb{R}^2 \ = \mathcal{G} \times_H \mathbb{R}^{n+2}/\mathbb{R}. \end{cases} \tag{23}$$

Using the $H$-module isomorphisms $(\mathfrak{q}^E_{-1} \oplus \mathfrak{q})/\mathfrak{h} \cong \mathbb{R}^{2*}$ and $(\mathfrak{q}^V_{-1} \oplus \mathfrak{q})/\mathfrak{h} \cong \mathbb{R}^{n+2}/\mathbb{R}$, we obtain (21) as desired. This proves the first assertion of the theorem.

Assume now that $(M, E, V)$ arises from a path geometry $(M', E', V')$ without auxiliary orientation where $M \to M'$ is the two-fold covering encoding the pair of orientations on $E'^* \otimes \wedge^n V'$, $E$ is the pullback of $E'$ and $V$ is the pullback of $V'$. Applying (20) to $(M, E, V)$ results in $(\tilde{M}, \tilde{E}, \tilde{F})$ where $\tilde{M} = GL_+(E^* \otimes \wedge^n V) = GL(E'^* \otimes \wedge^n V')$. It suffices to construct an almost Grassmann structure $(\tilde{M}', \tilde{E}', \tilde{F}')$ where $\tilde{M}' = GL(E'^* \otimes \wedge^n V')/\mathbb{Z}_2$ and $(\tilde{M}, \tilde{E}, \tilde{F})$ is the pullback of $(\tilde{M}', \tilde{E}', \tilde{F}')$.

Indeed, denote by $\hat{Q} \subseteq PGL(n+2, \mathbb{R})$ the maximal subgroup with Lie algebra $\mathfrak{q}$. We lift $\hat{Q}$ along the natural projection $GL(n+2, \mathbb{R}) \twoheadrightarrow PGL(n+2, \mathbb{R})$ to the block upper triangular subgroup of $GL(n+2, \mathbb{R})$ with block diagonal entries $\mathrm{diag}(a, b, C)$ such that $b > 0$ and $a\, b\, \det C = \pm 1$. This choice provides a unique representative for each element of $\hat{Q}$. This way, we obtain an identification

$$\hat{Q} \cong \mathrm{diag}(\pm 1, 1, \mathbb{I}_n)\, Q = Q\, \mathrm{diag}(\pm 1, 1, \mathbb{I}_n) \subseteq GL(n+2, \mathbb{R}).$$

Now let $\hat{H} \trianglelefteq \hat{Q}$ be the kernel of

$$\hat{Q} \twoheadrightarrow GL\left(\left((\mathfrak{q}^E_{-1} \oplus \mathfrak{q})/\mathfrak{q}\right)^* \otimes \wedge^n \left((\mathfrak{q}^V_{-1} \oplus \mathfrak{q})/\mathfrak{q}\right)\right).$$

The identification above restricted to $\hat{H}$ induces

$$\hat{H} \cong \mathrm{diag}(\pm 1, 1, \mathbb{I}_n)\, H = H\, \mathrm{diag}(\pm 1, 1, \mathbb{I}_n).$$

Recall from, e.g., [11, Proposition 4.4.3] that $(M', E', V')$ comes with a unique normal regular parabolic geometry $(\mathcal{G}' \to M', \omega')$ of type $(PGL(n+2, \mathbb{R}), \hat{Q})$. Then

$$(\mathcal{G}' \to \mathcal{G}'/Q, \omega')$$

is the normal regular parabolic geometry of type $(G, Q)$ corresponding to $(M, E, V)$. Similar to the previous step, the Cartan geometry

$$(\mathcal{G}' \to \mathcal{G}'/\hat{H}, \omega')$$

has the desired underlying almost Grassmann structure $(\tilde{M}', \tilde{E}', \tilde{F}')$ which clearly pulls back to the underlying geometry $(\tilde{M}, \tilde{E}, \tilde{F})$ of

$$(\mathcal{G}' \to \mathcal{G}'/H, \omega').$$

This completes the proof of the theorem. Note that the second half of the proof can be rephrased as an extension functor extending the Fefferman-type construction (19); see [13]. □

***Remark 12*** Although we have constructed $\tilde{E}^*$ and $\tilde{F}$ as distributions $(Tp)^{-1}(E)$ and $(Tp)^{-1}(V)$ of $\tilde{M}$, respectively, it offers more clarity to think of $\tilde{E}^*$ and $\tilde{F}$ as auxiliary vector bundles. Observe that for all $a \in \mathbb{R}$, $(0, a) \in \mathbb{R}^{2*}$ and $(a, 0^n)^t \in \mathbb{R}^{n+2}/\mathbb{R}$ are fixed by $H$. We specify the sections

$$\varphi = \mathcal{G} \times_H (0, 1) \in \Gamma(\tilde{E}^*) \quad \text{and} \quad \eta = \mathcal{G} \times_H (-1, 0^n)^t \in \Gamma(\tilde{F});$$

see (23). Denote by $\zeta_a \in \mathfrak{X}(\tilde{M})$ the fundamental vector field of $p : \tilde{M} \to M$ generated by $a \in \mathbb{R}$. Then the defining isomorphism $\tilde{E}^* \otimes \tilde{F} \cong T\tilde{M}$ of $(\tilde{M}, \tilde{E}, \tilde{F})$ restricts to

$$\tilde{E}^* \otimes \eta = (Tp)^{-1}(E), \quad \varphi \otimes \tilde{F} = (Tp)^{-1}(V) \quad \text{and} \quad \varphi \otimes \eta = \zeta_{-1}.$$

Let $N$ be an $(n+1)$-dimensional manifold with full frame bundle $\mathcal{P}N$. For $s \in \mathbb{R}$, recall that the $(-\frac{s}{n+2})$-density bundle

$$\mathcal{E}(s) = \mathcal{P}N \times_{\rho(s)} \mathbb{R}$$

is associated to the representation

$$\begin{aligned} \rho(s) : GL(n+1, \mathbb{R}) &\to GL_+(\mathbb{R}) \\ C &\mapsto |\det(C)|^{s/(n+2)}; \end{aligned}$$

see, e.g., [11, Subsection 4.1.5]. For example, if $N$ is an oriented projective structure, the standard tractor bundle $\mathcal{T}$ of the geometry comes with a canonical short exact sequence

$$0 \to \mathcal{E}(-1) \to \mathcal{T} \to TN \otimes \mathcal{E}(-1) \to 0;$$

see, e.g., [11, Subsection 5.2.6]. Now let $N$ be any $(n+1)$-dimensional manifold. Denote by $\mathcal{P}TN \to N$ the projectivized tangent bundle of $N$ and

$$V \subseteq \mathcal{H} \subseteq T\mathcal{P}TN$$

the vertical distribution and the tautological distribution, respectively. Recall that any line subbundle $E \subseteq \mathcal{H}$ such that $\mathcal{H} = E \oplus V$ defines a full path geometry $(\mathcal{P}TN, E, V)$.

In [14, Section 3], the authors associate to every full path geometry $(\mathcal{P}TN, E, V)$ with oriented $N$ a Segre structure on the weighted slit tangent bundle $TN \otimes \mathcal{E}(-1) \setminus \{0\}$. We show that this space can be identified with our Fefferman-type space. Assume now that $N$ is not necessarily oriented. We denote by

$$\widehat{\mathcal{P}TN} = (\wedge^{n+1} \mathcal{H} \setminus \{0\})/\mathbb{R}_+ \to \mathcal{P}TN$$

the two-fold covering encoding the pair of orientations on $\mathcal{H} \subseteq T\mathcal{P}TN$. Then every distribution of $\mathcal{P}TN$ canonically lifts to a distribution of $\widehat{\mathcal{P}TN}$, which we denote by the same notation. In particular, every full path geometry $(\mathcal{P}TN, E, V)$ canonically lifts to the path geometry $(\widehat{\mathcal{P}TN}, E, V)$ where

$$E \oplus V = \mathcal{H} \subseteq T\widehat{\mathcal{P}TN}$$

is oriented.

**Theorem 13** *Given a full path geometry* $(\mathcal{P}TN, E, V)$*, the Fefferman-type space of* $(\widehat{\mathcal{P}TN}, E, V)$ *is isomorphic to* $TN \otimes \mathcal{E}(-1) \setminus \{0\}$*, viewed as a fiber bundle over* $N$ *with standard fiber* $\mathbb{R}^{n+1} \setminus \{0\}$.

***Proof*** We first prove the statement of the theorem for every full path geometry on $N$ that arises from a certain correspondence space. Indeed, up to shrinking $N$, we may assume without loss of generality that there exists a normal parabolic geometry

$$(\mathcal{G} \to N, \omega)$$

of type $(SL(n+2, \mathbb{R}), R)$. Here, $R \subseteq SL(n+2, \mathbb{R}) = G$ is the stabilizer of $(\mathbb{R}, 0^{n+1})^t \subseteq \mathbb{R}^{n+2}$. We denote the Lie algebra of $R$ by $\mathfrak{r}$. It is well-known that $\mathcal{E}(s) \cong \mathcal{G} \times_R \mathbb{R}$ where the representation of $R$ on $\mathbb{R}$ is given by

$$R \ni \left(\begin{array}{c|c|c} a & * & * \\ \hline 0 & * & * \\ \hline & & \\ 0 & * & * \\ & & \end{array}\right) \mapsto |a|^{-s} \in GL_+(\mathbb{R}).$$

The block size is $(1, 1, n) \times (1, 1, n)$. Note that the above homomorphism is the result of the composition

$$R \xrightarrow{Ad} GL(\mathfrak{g}/\mathfrak{r}) = GL(n+1, \mathbb{R}) \xrightarrow{\rho(s)} GL_+(\mathbb{R}).$$

On the other hand, consider the representation of $R$ on

$$\mathbb{W} = (\mathfrak{g}/\mathfrak{r}, \underline{Ad}) \otimes (\mathbb{R}, \rho(-1)).$$

To gain intuition, note that $\mathbb{W} \cong \mathbb{R}^{n+2}/\mathbb{R}$ as $R$-modules. Observe that $H$ is the isotropic subgroup of some nonzero base point whose $R$-orbit contains all nonzero elements in the module. Hence $R/H \cong \mathbb{W} \setminus \{0\}$ and so

$$\mathcal{G}/H = \mathcal{G} \times_R R/H = \mathcal{G} \times_R \mathbb{W} \setminus \{0\} \cong TN \otimes \mathcal{E}(-1) \setminus \{0\}$$

as fiber bundles over $N$.

Analogous to, e.g., [11, Theorem 4.4.3], the representation of $R$ on $\mathfrak{g}/\mathfrak{r}$ gives rise to the identification

$$\mathcal{G}/Q = \mathcal{G} \times_R R/Q = \widehat{\mathcal{P}TN}$$

of the correspondence space. Moreover,

$$(\mathcal{G} \to \mathcal{G}/Q, \omega)$$

is the normal regular parabolic geometry canonically associated to the canonical lifting of a full path geometry $(\mathcal{P}TN, E, V)$ to $\widehat{\mathcal{P}TN}$. Its Fefferman-type space $\mathcal{G}/H$ admits the identification as desired.

It has been proven in Theorem 11 that given any full path geometry $(\mathcal{P}TN, E, V)$, the Fefferman-type space of $(\widehat{\mathcal{P}TN}, E, V)$ is isomorphic to the principal $\mathbb{R}_+$-bundle $GL_+((\mathcal{H}/V)^* \otimes \wedge^n V)$ over $\widehat{\mathcal{P}TN}$. We conclude from the above that $GL_+((\mathcal{H}/V)^* \otimes \wedge^n V)$, viewed as a fiber bundle over $N$, is isomorphic to $TN \otimes \mathcal{E}(-1) \setminus \{0\}$. This proves the statement of the theorem. □

We need the following technical result to study the construction using Cartan geometry.

**Lemma 14** *Let $(\mathcal{G} \to M, \omega)$ be a normal regular parabolic geometry of type $(G, Q)$ that is canonically associated to a path geometry $(M, E, V)$ where $M$ is $(2n+1)$-dimensional and $E \oplus V$ is oriented. Consider the $Q$-submodules*

$$A = \left\{ \left( \begin{array}{c|c|c} 0 & 0 & * \\ \hline 0 & 0 & * \\ \hline & & \\ 0 & 0 & * \\ & & \end{array} \right) \right\} \subseteq B = \left\{ \left( \begin{array}{c|c|c} 0 & * & * \\ \hline 0 & * & * \\ \hline & & \\ 0 & * & * \\ & & \end{array} \right) \right\} \subseteq \mathfrak{g} = \mathfrak{sl}(n+2, \mathbb{R})$$

*where the block size is $(1, 1, n) \times (1, 1, n)$. The following properties hold for the curvature $\kappa$ of the Cartan connection $\omega$.*

*(i) $\kappa \in \Omega^2(M, \mathcal{G} \times_Q B)$, $\kappa(V, TM) \subseteq \Gamma(\mathcal{G} \times_Q A)$ and $\kappa(E \oplus V, E \oplus V) = 0$.*
*(ii) If $(M, E, V)$ is torsion-free, then in addition to (i), we have that $\kappa \in \Omega^2(M, \mathcal{G} \times_Q A)$ and $\kappa(E, TM) = 0$.*

***Proof*** We use the improved Bianchi identity given in Lemma 5. Recall in Convention 6 the $|2|$-grading decomposition of $\mathfrak{g}$ and from [11, Subsection 4.4.3] that the harmonic curvature

$$\kappa_H : \mathcal{G} \to \ker \partial^* / \operatorname{im} \partial^* \cong \ker \Box \subseteq \wedge^2 \mathfrak{q}_+ \otimes \mathfrak{g}$$

is given by $\kappa_H = (\tau, \rho)$ where

$$\tau : \mathcal{G} \to (\mathfrak{q}_1^E \wedge \mathfrak{q}_2 \otimes \mathfrak{q}_{-1}^V) \cap \ker \Box \subseteq (\mathfrak{q}_1^E \oplus \mathfrak{q}_2) \wedge \mathfrak{q}_2 \otimes B$$

and

$$\rho : \mathcal{G} \to (\mathfrak{q}_1^V \wedge \mathfrak{q}_2 \otimes \mathfrak{q}_0) \cap \ker \Box \subseteq \mathfrak{q}_1^V \wedge \mathfrak{q}_2 \otimes A.$$

Here, we have used that

$$(\mathfrak{q}_1^V \wedge \mathfrak{q}_2 \otimes \mathfrak{q}_0) \cap \ker \Box \subseteq (\mathfrak{q}_1^V \wedge \mathfrak{q}_2 \otimes \mathfrak{q}_0) \cap \ker \partial^* \subseteq \mathfrak{q}_1^V \wedge \mathfrak{q}_2 \otimes \mathfrak{q}_0^{ss},$$

where $\mathfrak{q}_0^{ss} = [\mathfrak{q}_0, \mathfrak{q}_0] \subseteq A$.

Consider the $Q$-module

$$\mathbb{F} = (\mathfrak{q}_1^V \wedge \mathfrak{q}_2 \otimes A) \oplus (\mathfrak{q}_1^E \wedge \mathfrak{q}_2 \otimes B) \oplus (\mathfrak{q}_2 \wedge \mathfrak{q}_2 \otimes B).$$

Then the images of $\kappa_H$ lie in $\mathbb{F}$. We observe that for $\varphi, \psi \in \mathbb{F}$, $\iota_\varphi \psi \in (\mathfrak{q}_1^E \oplus \mathfrak{q}_2) \wedge \mathfrak{q}_2 \wedge \mathfrak{q}_2 \otimes A$. Applying $\partial^*$ on both sides we get $\partial^*(\iota_\varphi \psi) \in (\mathfrak{q}_1^E \oplus \mathfrak{q}_2) \wedge \mathfrak{q}_2 \otimes \mathfrak{q}_2 \subseteq \mathbb{F}$. Hence $\mathbb{F}$ is stable under $\mathbb{F}$-insertions and so $\kappa$ is a section of $\Gamma(\mathcal{G} \times_Q \mathbb{F})$, proving the claim in (i).

Now suppose that $(M, E, V)$ is torsion-free, i.e., $\kappa_H = \rho$. Then the images of $\kappa_H$ lie in the $Q$-module

$$\mathbb{E} = (\mathfrak{q}_1^V \oplus \mathfrak{q}_2) \wedge \mathfrak{q}_2 \otimes A.$$

Observe that $\mathbb{E}$ is trivially stable under $\mathbb{E}$-insertions because for every $\varphi, \psi \in \mathbb{E}$, $\iota_\varphi \psi = 0$. Hence $\kappa$ is a section of $\Gamma(\mathcal{G} \times_Q \mathbb{E})$, proving the claim in (ii). □

**Proposition 15** *Notation as at the beginning of the section. Then the Cartan connection $\tilde{\omega}$ is automatically normal. That is, the Fefferman-type construction from path geometries to almost Grassmann structures is normal.*

***Proof*** Denote by $\partial^*$ and $\tilde{\partial}^*$ the Kostant codifferentials of $(\mathfrak{g}, \mathfrak{q})$ and of $(\tilde{\mathfrak{g}}, \tilde{\mathfrak{p}})$, respectively. Consider

$$\phi \in \wedge^2 \mathfrak{q}_+ \otimes \mathfrak{g} \quad \text{and} \quad \tilde{\phi} \in \wedge^2 \tilde{\mathfrak{p}}_+ \otimes \tilde{\mathfrak{g}}$$

such that

$$\tilde{\phi} = i' \circ \phi \circ \wedge^2 \pi;$$

where $\pi$ is as at the beginning of the section.

Let $\alpha : \mathfrak{g} \to \tilde{\mathfrak{g}}$ be the embedding of $\mathfrak{g} = \mathfrak{sl}(n+2, \mathbb{R})$ to the Lie subalgebra of $\tilde{\mathfrak{g}} = \mathfrak{sl}(n+3, \mathbb{R})$ with zeros in the third row and the third column and $B \subseteq \mathfrak{g}$ be the $Q$-submodule defined in Lemma 14. We show that if $\phi \in \wedge^2 \mathfrak{q}_+ \otimes B$ and $\phi(\mathfrak{q}_{-1}^E, \mathfrak{q}_{-1}^V) = 0$, then

$$\tilde{\partial}^* \tilde{\phi} = \alpha \circ (\partial^* \phi) \circ \pi.$$

Indeed, observe that the dual of $\pi : \tilde{\mathfrak{g}}/\tilde{\mathfrak{p}} \twoheadrightarrow \mathfrak{g}/\mathfrak{q}$ is $\pi^* : \mathfrak{q}_+ \to \tilde{\mathfrak{p}}_+$ where, as matrices,

$$\pi^*(Z)_{ij} = Z_{i,j+1}, \quad Z \in \mathfrak{q}_+. \tag{24}$$

That is, we have $\langle Z, \pi(\tilde{X}) \rangle = \langle \pi^*(Z), \tilde{X} \rangle$ for all $Z \in \mathfrak{q}_+$ and $\tilde{X} \in \tilde{\mathfrak{g}}/\tilde{\mathfrak{p}}$. Observe that for all $Z \in \mathfrak{q}_+$ and $W \in B$, there holds

$$[\pi^*(Z), i'(W)] = \alpha([Z, W]). \tag{25}$$

Moreover, let $X \in \mathfrak{g}$ be any element and $\{X^i\}$ a basis of $\mathfrak{g}/\mathfrak{q}$ with dual basis $\{Z_i\}$ of $\mathfrak{q}_+$. Then

$$(\partial^* \phi)(X) = 2 \sum_i [Z_i, \phi(X, X^i)] - \sum_i \phi([Z_i, X], X^i) = 2 \sum_i [Z_i, \phi(X, X^i)];$$

see (8). Note that the assumption $\phi(\mathfrak{q}^E_{-1}, \mathfrak{q}^V_{-1}) = 0$ implies that $\phi([Z_i, X], X^i) = 0$ for all pairs $(X^i, Z_i)$.

Let $(\{\tilde{X}^i\}, \tilde{Y})$ be a basis of $\tilde{\mathfrak{g}}/\tilde{\mathfrak{p}}$ such that $\pi(\tilde{X}^i) = X^i$ for all $i$ and $\tilde{Y} \in \ker \pi$. Then its dual basis is $(\{\pi^* Z_i\}, \tilde{W}) \in \tilde{\mathfrak{p}}_+$ where $\tilde{W}$ is uniquely determined. For any $\tilde{X} \in \tilde{\mathfrak{g}}/\tilde{\mathfrak{p}}$, write $X = \pi(\tilde{X})$. Then, using (8), the assumption $\pi(\tilde{Y}) = 0$ and (25),

$$\begin{aligned}(\tilde{\partial}^* \tilde{\phi})(\tilde{X}) &= 2\left([\tilde{W}, i' \circ \phi(X, 0)] + \sum_i [\pi^*(Z_i), i' \circ \phi(X, X^i)]\right)\\ &= 2\sum_i \alpha([Z_i, \phi(X, X^i)])\\ &= (\partial^* \phi)(X).\end{aligned}$$

This proves the claim.

Using Lemma 14 (i) and the above claim, we conclude that

$$(\tilde{\partial}^* \tilde{\kappa})(j(u)) = \alpha \circ (\partial^* \kappa)(u) \circ \pi = 0$$

for every $u \in \mathcal{G}$; where $j$, $\kappa$ and $\tilde{\kappa}$ are as at the beginning of the section. This completes the proof of the proposition. □

Consequently, the natural map

$$j : \mathcal{G} \to \tilde{\mathcal{G}} = \mathcal{G} \times_H \tilde{P}$$

defining $\tilde{\mathcal{G}}$ relates Cartan objects of $(M, E, V)$ to those of $(\tilde{M}, \tilde{E}, \tilde{F})$ in a canonical way. For example, the two equalities at the beginning of Sect. 3 are relations between the normal regular Cartan connection $\omega$ of the initial path geometry and the normal Cartan connection $\tilde{\omega}$ of the resulting almost Grassmann structure and between the curvature $\kappa$ of $\omega$ and the curvature $\tilde{\kappa}$ of $\tilde{\omega}$, respectively. A tractor bundle on $(\tilde{M}, \tilde{E}, \tilde{F})$ associated to a representation of $\tilde{G}$ is the pullback of the tractor bundle on $(M, E, V)$ associated to the representation induced by $i : G \to \tilde{G}$. Moreover, recall that the canonical curves on $(M, E, V)$ are projections of the flows $Fl^{\omega^{-1}(X)}(u)$ along $\mathcal{G} \to M$ where $X \in \mathfrak{g}$ and $u \in \mathcal{G}$. Observe that

$$j \circ Fl^{\omega^{-1}(X)}(u) = Fl^{\tilde{\omega}^{-1}(i(X))}(j(u)).$$

It follows that every canonical curve of $(M, E, V)$ is lifted along $\tilde{M} \twoheadrightarrow M$ to a canonical curve of $(\tilde{M}, \tilde{E}, \tilde{F})$.

Additionally, the following properties hold for the Fefferman-type construction we specified. Notation as at the beginning of the paragraph. Then

$$\tilde{\mathcal{T}} = \tilde{\mathcal{G}} \times_{\tilde{P}} \mathbb{R}^{n+3} \quad \text{and} \quad \tilde{\mathcal{G}} \times_{\tilde{P}} \tilde{\mathfrak{g}} = \mathfrak{sl}(\tilde{\mathcal{T}})$$

are the standard tractor bundle and the adjoint tractor bundle of $(\tilde{M}, \tilde{E}, \tilde{F})$, respectively. Then

$$\tilde{\kappa} \in \Omega^2(\tilde{M}, \mathfrak{sl}(\tilde{\mathcal{T}})).$$

Recall the natural short exact sequences

$$0 \to \tilde{E} \to \tilde{\mathcal{T}} \to \tilde{F} \to 0 \quad \text{and} \quad 0 \to \tilde{F}^* \to \tilde{\mathcal{T}}^* \to \tilde{E}^* \to 0.$$

They follow easily from Lemma 7. Consider the nowhere vanishing sections

$$\varphi = \mathcal{G} \times_H (0, 1) \in \Gamma(\tilde{E}^*) \quad \text{and} \quad \eta = \mathcal{G} \times_H (-1, 0^n)^t \in \Gamma(\tilde{F})$$

defined in Remark 12. Moreover, denote the torsion of $(\tilde{M}, \tilde{E}, \tilde{F})$ by $\tilde{\tau}$. The contraction of $\tilde{\tau} \otimes \tilde{\tau}$ specified in Lemma 8 results in

$$tr(\iota_{\tilde{\tau}}\tilde{\tau}) \in \Gamma(\mathrm{Sym}^2\, T^*\tilde{M}).$$

**Theorem 16** *Consider the Fefferman-type construction*

$$(M, E, V) \mapsto (\tilde{M}, \tilde{E}, \tilde{F})$$

*of almost Grassmann structures from path geometries. The following properties hold.*

*(i)* $tr(\iota_{\tilde{\tau}}\tilde{\tau}) = 0$.
*(ii)* $(\tilde{M}, \tilde{E}, \tilde{F})$ *is torsion-free if and only if* $(M, E, V)$ *is torsion-free.*
*(iii) There are parallel sections* $\mu \in \Gamma(\tilde{\mathcal{T}})$ *and* $\nu \in \Gamma(\tilde{\mathcal{T}}^*)$ *with projections* $\eta \in \Gamma(\tilde{F})$ *and* $\varphi \in \Gamma(\tilde{E}^*)$*, respectively. Moreover, there holds* $\mu(\nu) = 1$.
*(iv) Denote by*

$$\ker \varphi \subseteq \tilde{E} \subseteq \mathrm{span}(\eta) + \tilde{E} \subseteq \tilde{\mathcal{T}}$$

*the natural filtration of ranks* $1 \leq 2 \leq 3 \leq n+3$ *and* $\tilde{\mathcal{D}} = \mathrm{span}(\varphi \otimes \tilde{F}, \tilde{E}^* \otimes \eta) \subseteq T\tilde{M}$. *Then* $\tilde{\kappa}$ *has the following properties.*

$$\begin{cases} \tilde{\kappa}(\varphi \otimes \eta, T\tilde{M}) & = 0 \\ \tilde{\kappa}(\tilde{\mathcal{D}}, \tilde{\mathcal{D}}) & = 0 \\ \tilde{\kappa}(\varphi \otimes \tilde{F}, T\tilde{M}) & \subseteq \{\phi \in \Gamma(\mathfrak{sl}(\tilde{\mathcal{T}})) : \phi(\mathrm{span}(\eta) + \tilde{E}) = 0\} \quad \textit{and} \\ \tilde{\kappa}(T\tilde{M}, T\tilde{M}) & \subseteq \{\phi \in \Gamma(\mathfrak{sl}(\tilde{\mathcal{T}})) : \phi(\ker \varphi) = 0\}. \end{cases}$$

*If* $(M, E, V)$ *is torsion-free, then in addition to the above we have that*

$$\begin{cases} \tilde{\kappa}(\tilde{E}^* \otimes \eta, T\tilde{M}) & = 0 \quad \textit{and} \\ \tilde{\kappa}(T\tilde{M}, T\tilde{M}) & \subseteq \{\phi \in \Gamma(\mathfrak{sl}(\tilde{\mathcal{T}})) : \phi(\mathrm{span}(\eta) + \tilde{E}) = 0\}. \end{cases}$$

***Proof*** Notation as in the first paragraph of the proof of Proposition 15. Recall that $\tilde{\kappa}$ is completely determined by the curvature $\kappa$ of $(\mathcal{G} \to M, \omega)$ via the relation

$$\tilde{\kappa}(j(u)) = i' \circ \kappa(u) \circ \wedge^2 \pi$$

for all $u \in \mathcal{G}$. Thus the properties of $\kappa$ given in Lemma 14 (i) give rise to the properties of $\tilde{\kappa}$ in (iv) in the statement of the theorem. Recall from Remark 12 that $\varphi \otimes \eta$ spans the vertical bundle of $p : \tilde{M} \to M$.

Let $A \subseteq \mathfrak{g}$ be the $Q$-submodule defined in Lemma 14. That $i'^{-1}(\tilde{\mathfrak{p}}) = \mathfrak{h} \subseteq \mathfrak{q}$ immediately implies that if $\tilde{\phi} \in \wedge^2 \tilde{\mathfrak{p}}_+ \otimes \tilde{\mathfrak{p}}$, then $\phi \in \wedge^2 \mathfrak{q}_+ \otimes \mathfrak{q}$. Moreover, that $i'(A) \subseteq \tilde{\mathfrak{p}}$ immediately implies that if $\phi \in \wedge^2 \mathfrak{q}_+ \otimes A$, then $\tilde{\phi} \in \wedge^2 \tilde{\mathfrak{p}}_+ \otimes \tilde{\mathfrak{p}}$. Using Lemma 14 and these two relations, we conclude that (ii) holds.

Observe that Theorem 4 applies to our Fefferman-type functor from Cartan geometries of type $(G, Q)$ to Cartan geometries of type $(\tilde{G}, \tilde{P})$ defined by $i : G \to \tilde{G}$ where we choose the tractor bundle $\tilde{\mathcal{W}} = \tilde{\mathcal{T}} \oplus \tilde{\mathcal{T}}^*$ and the fiber subbundle $\tilde{\mathcal{O}} = \tilde{\mathcal{G}} \times_{\tilde{P}} (\tilde{P}\, \tilde{w}_0)$ determined by the point

$$\tilde{w}_0 = \big((0, 0, -1, 0^n)^t, (0, 1, -1, 0^n)\big) \in \mathbb{R}^{n+3} \oplus \mathbb{R}^{(n+3)*}$$

fixed by $H$. As explained in the proof of Theorem 4, the constant map $\mathcal{G} \to \{\tilde{w}_0\}$ uniquely extends to a pair

$$(\mu, \nu) \in C^\infty(\tilde{\mathcal{G}}, \mathbb{R}^{n+3} \oplus \mathbb{R}^{(n+3)*})^{\tilde{P}} = \Gamma(\tilde{\mathcal{T}} \oplus \tilde{\mathcal{T}}^*)$$

of parallel sections. Clearly, $(\mu, \nu)$ projects to $(\eta, \varphi)$ and $\mu(\nu) = 1$. This proves (iii).

It remains to prove (i). For every $\tilde{x} \in \tilde{M}$, we may choose $u \in \mathcal{G}$ such that $j(u) \in \tilde{\mathcal{G}}$ projects to $\tilde{x}$. Let $\tilde{\sigma} : \tilde{\mathcal{G}}/\tilde{P}_+ \to \tilde{\mathcal{G}}$ be a Weyl structure such that $j(u)$ lies in its image. Consider the component

$$\tilde{W} \in \Omega^2(\tilde{M}, L(\tilde{E}, \tilde{E}))$$

of the Weyl tensor of $\tilde{\sigma}$; see (14). It follows from the fourth property of $\tilde{\kappa}$ in (iv) that

$$\tilde{W}_{\tilde{x}}(T_{\tilde{x}}\tilde{M}, T_{\tilde{x}}\tilde{M}) \subseteq \varphi(\tilde{x}) \otimes \tilde{E}_{\tilde{x}}$$

and from the third property of $\tilde{\kappa}$ in (iv) that

$$\tilde{W}_{\tilde{x}}(\varphi(\tilde{x}) \otimes \tilde{F}_{\tilde{x}}, T_{\tilde{x}}\tilde{M}) = 0.$$

Recall that in Lemma 8 we fix a Ricci-type contraction $tr(\tilde{W})$ of $\tilde{W}$ and provide a relation between $tr(\tilde{W})$ and $tr(\iota_{\tilde{\tau}}\tilde{\tau})(\tilde{x})$, which implies that $tr(\iota_{\tilde{\tau}}\tilde{\tau})(\tilde{x}) = 0$ if and only if $tr(\tilde{W}_{\tilde{x}}) = 0$. In fact, the above two properties on $\tilde{W}_{\tilde{x}}$ together imply that $tr(\tilde{W}_{\tilde{x}}) = 0$. This proves (i). □

From (iv) in Theorem 16 we conclude that

$$\tilde{\tau} \in \Gamma\big(L\big(\wedge^2 (T\tilde{M}/(\varphi \otimes \tilde{F})), \varphi \otimes \tilde{F}\big)\big). \tag{26}$$

Observe that $\tilde{\tau}$ is completely determined by the projection of the regular normal Cartan curvature $\kappa \in \Omega^2(M, \mathcal{G} \times_Q \mathfrak{g})$ of the path geometry $(M, E, V)$ in $\Omega^2(M, \mathcal{G} \times_Q \mathfrak{g}/\mathfrak{h})$;

see Convention 6. If $(M, E, V)$ is torsion-free, then $\tilde{\tau} = 0$ by Theorem 16 (ii) and the harmonic curvature $\tilde{\rho} \in \Omega^2(\tilde{M}, \mathfrak{sl}(\tilde{F}))$ of $(\tilde{M}, \tilde{E}, \tilde{F})$ equals the projection of $\tilde{\kappa} \in \Omega^2(\tilde{M}, \tilde{\mathcal{G}} \times_{\tilde{P}} \tilde{\mathfrak{p}})$ to $\Omega^2(\tilde{M}, \tilde{\mathcal{G}} \times_{\tilde{P}} \tilde{\mathfrak{p}}/\tilde{\mathfrak{p}}_+)$; see Convention 6 and, e.g., [11, Theorem 3.1.12]. In particular,

$$\tilde{\rho} \in \Gamma\big(L\big(\wedge^2 (T\tilde{M}/(\tilde{E}^* \otimes \eta)), \mathfrak{sl}(\tilde{F})\big)\big) \tag{27}$$

and is completely determined by the projection of $\kappa \in \Omega^2(M, \mathcal{G} \times_Q \mathfrak{q})$ to $\Omega^2(M, \mathcal{G} \times_Q \mathfrak{q}/\mathfrak{q}_2)$.

Combining (i) and (ii) in Theorem 16, we see that there are plenty of almost Grassmann structures with non-vanishing torsion $\tau$ such that $tr(\iota_\tau \tau) = 0$. Thus the property $tr(\iota_\tau \tau) = 0$ is a nontrivial weakening of torsion-freeness.

The tensor field $tr(\iota_\tau \tau)$ turns out to be important to almost Grassmann structures: We prove in Proposition 24 below that a Fefferman-type construction from an almost Grassmann structure of type $(2, n)$ to an almost Grassmann structure of type $(2, n+1)$ is normal if and only if the initial almost Grassmann structure of type $(2, n)$ has the property that $tr(\iota_\tau \tau) = 0$. It is proven in [18, Theorem 8 and Proposition 11] that parallel standard tractors (respectively, parallel standard cotractors) of an almost Grassmann structure $(M, E, F)$ are in one-to-one correspondence with sections of $F$ (respectively, of $E^*$) which are parallel for some Weyl connection and annihilate $tr(\iota_\tau \tau)$.

Conversely, we want to characterize all almost Grassmann structures $(\tilde{M}, \tilde{E}, \tilde{F})$ that are locally the result of the Fefferman-type construction. That is, there is an open cover of $\tilde{M}$ such that the restriction of $(\tilde{M}, \tilde{E}, \tilde{F})$ to every open set in the cover is isomorphic to an open subspace of an almost Grassmann structure resulted from the Fefferman-type construction. Here, the restriction of the almost Grassmann structure

$$(\tilde{M}, \tilde{E}, \tilde{F}) = (\tilde{M}, p : \tilde{E} \to \tilde{M}, q : \tilde{F} \to \tilde{M})$$

to an open subset $\tilde{U} \subseteq \tilde{M}$ refers to the almost Grassmann structure $(\tilde{U}, p^{-1}(\tilde{U}), q^{-1}(\tilde{U}))$, which is said to be an open subspace of $(\tilde{M}, \tilde{E}, \tilde{F})$.

Let $(\tilde{M}, \tilde{E}, \tilde{F})$ be an almost Grassmann structure of type $(2, n+1)$ with canonically associated normal parabolic geometry $(\tilde{\mathcal{G}} \to \tilde{M}, \tilde{\omega})$ of type $(\tilde{G}, \tilde{P})$. Denote by $\tilde{\kappa} \in \Omega^2(\tilde{M}, \tilde{\mathcal{G}} \times_{\tilde{P}} \tilde{\mathfrak{g}})$ the curvature of $\tilde{\omega}$ and by $\tilde{\mathcal{T}} = \tilde{\mathcal{G}} \times_{\tilde{P}} \mathbb{R}^{n+3}$ the standard tractor bundle of $(\tilde{M}, \tilde{E}, \tilde{F})$. Then $\tilde{\mathcal{T}}^*$ is the standard cotractor bundle of $(\tilde{M}, \tilde{E}, \tilde{F})$. Recall the natural projections $\tilde{\mathcal{T}} \twoheadrightarrow \tilde{F}$ and $\tilde{\mathcal{T}}^* \twoheadrightarrow \tilde{E}^*$. They may be obtained using Lemma 7 in a straightforward way.

**Theorem 17** *Notation as in the preceding paragraph. $(\tilde{M}, \tilde{E}, \tilde{F})$ is locally the result of the Fefferman-type construction on a path geometry if and only if the following hold.*

*(i) There are parallel sections $\mu \in \Gamma(\tilde{\mathcal{T}})$ and $\nu \in \Gamma(\tilde{\mathcal{T}}^*)$ with $\mu(\nu) = 1$. Moreover, the projection $\eta \in \Gamma(\tilde{F})$ of $\mu$ and $\varphi \in \Gamma(\tilde{E}^*)$ of $\nu$ are both nowhere vanishing.*

*(ii) $\tilde{\kappa}(\varphi \otimes \eta, \tilde{E}^* \otimes \tilde{F}) = 0$, $\tilde{\kappa}(\varphi \otimes \tilde{F}, \tilde{E}^* \otimes \eta) \subseteq \tilde{\mathcal{G}} \times_{\tilde{P}} \tilde{\mathfrak{p}}_+$ and $\varphi \otimes \tilde{F} \subseteq T\tilde{M}$ is an involutive distribution.*

*Moreover, if $(\tilde{M}, \tilde{E}, \tilde{F})$ is torsion-free, then (i) implies (ii) already.*

***Proof*** That all almost Grassmann structures resulting from our Fefferman-type construction satisfy (i) and (ii) has been proven in Theorem 16. We prove that the converse direction holds.

As mentioned in the proof of Theorem 16, Theorem 4 applies to our Fefferman-type functor from Cartan geometries of type $(G, Q)$ to Cartan geometries of type $(\tilde{G}, \tilde{P})$ defined by $i : G \to \tilde{G}$ where we choose the tractor bundle $\tilde{\mathcal{W}} = \tilde{\mathcal{T}} \oplus \tilde{\mathcal{T}}^*$ and the fiber subbundle $\tilde{\mathcal{O}} = \tilde{\mathcal{G}} \times_{\tilde{P}} (\tilde{P}\,\tilde{w}_0)$ determined by the point

$$\tilde{w}_0 = \big((0, 0, -1, 0^n)^t, (0, 1, -1, 0^n)\big) \in \mathbb{R}^{n+3} \oplus \mathbb{R}^{(n+3)*}$$

fixed by $H$. One could easily check that the sections of $\tilde{\mathcal{O}}$ are precisely pairs of sections $\mu' \in \Gamma(\tilde{\mathcal{T}})$ and $\nu' \in \Gamma(\tilde{\mathcal{T}}^*)$ such that $\nu'(\mu') = 1$ and that the projections $\eta' \in \Gamma(\tilde{F})$ and $\varphi' \in \Gamma(\tilde{E}^*)$ are both nowhere vanishing. Moreover, the distribution $\tilde{\mathcal{G}} \times_{\tilde{P}} (i'(\mathfrak{q}) + \tilde{\mathfrak{p}})/\tilde{\mathfrak{p}}$ is spanned by $\varphi' \otimes \eta'$.

Assume that (i) and (ii) hold. It follows from Theorem 4 that $(\tilde{\mathcal{G}} \to \tilde{M}, \tilde{\omega})$ is locally the result of the Fefferman-type functor applied to a parabolic geometry of type $(G, Q)$. Without loss of generality, assume that $(\tilde{\mathcal{G}} \to \tilde{M}, \tilde{\omega})$ is the result of the Fefferman-type functor applied to a parabolic geometry $(\mathcal{G} \to M, \omega)$ of type $(G, Q)$. Consider the distinguished distributions

$$E = \mathcal{G} \times_Q (\mathfrak{q}_{-1}^E \oplus \mathfrak{q})/\mathfrak{q} \quad \text{and} \quad V = \mathcal{G} \times_Q (\mathfrak{q}_{-1}^V \oplus \mathfrak{q})/\mathfrak{q}$$

on $M$. Denote by $p : \tilde{M} \to M$ the natural projection. All relations in Theorem 11 and Remark 12 hold. In particular, we have that

$$(Tp)^{-1}(E) = \tilde{E}^* \otimes \eta \quad \text{and} \quad (Tp)^{-1}(V) = \varphi \otimes \tilde{F}.$$

Thus the assumptions on $\tilde{\kappa}$ in (ii) imply that $(\mathcal{G} \to M, \omega)$ is a regular parabolic geometry with underlying path geometry $(M, E, V)$ with oriented $E \oplus V$.

Note that the characterizing properties in (i) and (ii) of the theorem is a subcollection of the properties listed in Theorem 16. This subcollection of conditions is sufficient to prove that $(\tilde{M}, \tilde{E}, \tilde{F})$ is the result of the Fefferman-type construction on $(M, E, V)$ but is insufficient to prove that $\omega$ is normal. In fact, let $\omega^{nor} \in \Omega^1(\mathcal{G}, \mathfrak{g})$ be the normal regular Cartan connection of type $(G, Q)$ canonically associated to the geometry $(M, E, V)$; see, e.g., [11, Proposition 3.1.13]. Recall that $\omega^{nor}$ is unique up to pullback by an automorphism on the principal bundle $\mathcal{G} \to M$; see, e.g., [11, Proposition 3.1.14]. We show that

- there holds $\omega^{nor} - \omega \in \Omega^1(\mathcal{G}, \mathfrak{q}^1)$ where $\mathfrak{q}^1 \subseteq \mathfrak{h} = i'^{-1}(\tilde{\mathfrak{p}})$; and
- if in addition $\tilde{\kappa}(\varphi \otimes \tilde{F}, \tilde{E}^* \otimes \eta) = 0$, we may choose $\omega^{nor} = \omega$.

Note that the assumption $\tilde{\kappa}(\varphi \otimes \tilde{F}, \tilde{E}^* \otimes \eta) = 0$ in the second claim above is also taken from the properties listed in Theorem 16. Indeed, we use the following notation from the proof of Proposition 15. Let $\partial^*$ and $\tilde{\partial}^*$ be the Kostant differentials of $(\mathfrak{g}, \mathfrak{q})$ and $(\tilde{\mathfrak{g}}, \tilde{\mathfrak{p}})$, respectively. Let $\pi^* : \mathfrak{q}_+ \to \tilde{\mathfrak{p}}_+$ be the dual map of $\pi : \tilde{\mathfrak{g}}/\tilde{\mathfrak{p}} \twoheadrightarrow \mathfrak{g}/\mathfrak{q}$, $\{X^i\} \subseteq \mathfrak{g}/\mathfrak{q}$ any basis with dual basis $\{Z_i\} \subseteq \mathfrak{q}_+$, $(\{\tilde{X}^i\}, \tilde{Y}) \subseteq \tilde{\mathfrak{g}}/\tilde{\mathfrak{p}}$ a basis such that $\pi(\tilde{X}_i) = X_i$ for all $i$ and $0 \neq \tilde{Y} \subseteq \ker \pi$ with dual basis $(\{\pi^*(Z_i)\}, \tilde{W}) \subseteq \tilde{\mathfrak{p}}_+$. Moreover, let

$$\beta : \tilde{\mathfrak{g}} \twoheadrightarrow \mathfrak{g}$$

be the linear map given by deleting the third row and summing up the second and the third column of $\tilde{\mathfrak{g}} = \mathfrak{sl}(n+3, \mathbb{R})$. Then we have the relation

$$\beta([\pi^*(Z), i'(W)]) = [Z, W]$$

for all $Z \in \mathfrak{q}_+$ and $W \in \mathfrak{g}$. Let $\phi \in \wedge^2 \mathfrak{q}_+ \otimes \mathfrak{g}$ and

$$\tilde{\phi} = i' \circ \phi \circ \wedge^2 \pi \in \wedge^2 \tilde{\mathfrak{p}}_+ \otimes \tilde{\mathfrak{g}}.$$

Recall that

$$(\partial^* \phi)(X) = 2 \sum_i [Z_i, \phi(X, X^i)] - \sum_i \phi([Z_i, X], X^i).$$

for all $X \in \mathfrak{g}$; see (8). Let $\tilde{X} = i'(X) \in \tilde{\mathfrak{g}}$. Since $\pi(\tilde{X} + \tilde{\mathfrak{p}}) = X + \mathfrak{q}$ and the pair $(\tilde{W}, \tilde{Y})$ of bases satisfies $[\tilde{W}, \tilde{\phi}(\tilde{X}, \tilde{Y})] = [\tilde{W}, i' \circ \phi(X, 0)] = 0$, we conclude that

$$2 \sum_i [Z_i, \phi(X, X^i)] = 2 \sum_i \beta([\pi^*(Z_i), \tilde{\phi}(\tilde{X}, \tilde{X}^i)]) = \beta \circ (\tilde{\partial}^* \tilde{\phi})(\tilde{X}). \tag{28}$$

On the other hand, we show that

$$-\sum_i \phi([Z_i, X], X^i) = \begin{cases} 0 & X \in \mathfrak{q}^{-1} \\ 2\,\phi(X_E, X_V) & X = [X_E, X_V], \quad X_E \in \mathfrak{q}^E_{-1}, X_V \in \mathfrak{q}^V_{-1}. \end{cases} \tag{29}$$

Indeed, we may assume without loss of generality that

$$Z_1 \in \mathfrak{q}^E_1, \quad Z_{1+j} \in \mathfrak{q}^V_1 \quad \text{and} \quad Z_{n+1+j} = [Z_1, Z_{1+j}] \in \mathfrak{q}_2$$

for all $1 \leq j \leq n$. Then

$$X^{n+1+j} = -[X^1, X^{1+j}].$$

If $X \in \mathfrak{q}^{-1}$, then $\phi([Z_i, X], X^i) \subseteq \phi(\mathfrak{q}, X^i) = 0$ for all $1 \leq i \leq 2n+1$; and if

$$X + \mathfrak{q} = X^{n+1+j_0} + \mathfrak{q} = -[X^1, X^{1+j_0}] + \mathfrak{q}$$

for some $1 \le j_0 \le n$, then

$$
\begin{aligned}
-\sum_i \phi([Z_i, [X^1, X^{1+j_0}]], X^i) &= \sum_i \phi([Z_i, X], X^i) \\
&= \phi([Z_1, X], X^1) + \phi([Z_{1+j_0}, X], X^{1+j_0}) \\
&= \phi(-X^{1+j_0}, X^1) + \phi(X^1, X^{1+j_0}) \\
&= 2\,\phi(X^1, X^{1+j_0}).
\end{aligned}
$$

This completes the proof of (29).

Now substitute $\phi = \kappa(u)$ and $\tilde{\phi} = \tilde{\kappa}(j(u))$ in the above where $\kappa$ is the curvature of $\omega$, $\tilde{\kappa}$ is the curvature of $\tilde{\omega}$, $j : \mathcal{G} \to \tilde{\mathcal{G}}$ the inclusion map and $u \in \mathcal{G}$ is any point. Since $\tilde{\kappa}(\tilde{E}^* \otimes \eta, \varphi \otimes \tilde{F}) \subseteq \tilde{\mathcal{G}} \times_{\tilde{P}} \tilde{\mathfrak{p}}_+$ and $i'^{-1}(\tilde{\mathfrak{p}}_+) = \mathfrak{q}_1^V \oplus \mathfrak{q}_2 \subseteq \mathfrak{q}^1$, we conclude that

$$
(\partial^*\kappa)(E \oplus V) = 0 \quad \text{and} \quad (\partial^*\kappa) \circ \mathcal{L}|_{E\times V} = 2\kappa|_{E\times V} \subseteq \mathcal{G} \times_Q \mathfrak{q}^1
$$

where $\mathcal{L} : \wedge^2(E \oplus V) \to \frac{TM}{E\oplus V}$ is the Levi bracket. In particular, $\partial^*\kappa$ has homogeneity $\ge 3$. It follows from, e.g., [11, Proposition 3.1.13] that $\omega^{nor} - \omega$ also has homogeneity $\ge 3$ and thus lies in $\Omega^1(\mathcal{G}, \mathfrak{q}^1)$. If in addition $\tilde{\kappa}(\tilde{E}^* \otimes \eta, \varphi \otimes \tilde{F}) = 0$, the above procedure implies that $\partial^*\kappa = 0$ and so $(\mathcal{G} \to M, \omega)$ is already the normal regular parabolic geometry canonically associated to $(M, E, V)$. This completes the proof of the above two claims.

Since $i'(\mathfrak{q}^1) \subseteq \tilde{\mathfrak{p}}$, the Fefferman-type functor applied to $(\mathcal{G} \to M, \omega^{nor})$ results in a parabolic geometry $(\tilde{\mathcal{G}} \to \tilde{M}, \widetilde{\omega^{nor}})$ such that $\widetilde{\omega^{nor}} - \tilde{\omega} \in \Omega^1(\tilde{\mathcal{G}}, \tilde{\mathfrak{p}})$. In particular, the underlying geometry of $(\tilde{\mathcal{G}} \to \tilde{M}, \widetilde{\omega^{nor}})$ is the same as the underlying geometry of $(\tilde{\mathcal{G}} \to \tilde{M}, \tilde{\omega})$; see, e.g., [11, Proposition 3.1.10 (1)]. Hence $(\tilde{M}, \tilde{E}, \tilde{F})$ is the result of the Fefferman-type construction on $(M, E, V)$. This proves the converse direction.

Finally, if $(\tilde{M}, \tilde{E}, \tilde{F})$ is torsion-free and satisfies (i), then $\tilde{\kappa}(\varphi \otimes \tilde{F}, \tilde{E}^* \otimes \eta)$ and $\tilde{\kappa}(\varphi\otimes\tilde{F}, \varphi\otimes\tilde{F})$ take values in $\tilde{\mathcal{G}}\times_{\tilde{P}}\tilde{\mathfrak{p}}_+$. In particular, $\varphi\otimes\tilde{F}$ is an involutive distribution. Moreover, observe that the trace-free part of $\mu\otimes\nu \in \mathfrak{gl}(\tilde{\mathcal{T}})$ is a parallel adjoint tractor. By [6, 3.5 Corollary], there holds $\tilde{\kappa}(\varphi \otimes \eta, \tilde{E}^* \otimes \tilde{F}) = 0$ and so we obtain (ii). This completes the prove of the theorem. □

Thus we have provided a full characterization of the Fefferman-type spaces we have introduced in a style analogous to the classical result in [8, Theorem 2.5] characterizing conformal structures that are locally constructed from CR structures by the existence of a parallel adjoint tractor. Recall from, e.g., [11, Subsection 4.4.5] that torsion-free path geometries are locally the correspondence space of a Grassmann structure. Hence given a Grassmann structure $(\tilde{M}, \tilde{E}, \tilde{F})$ of type $(2, n+1)$, there holds Theorem 17 (i) if and only if $(\tilde{M}, \tilde{E}, \tilde{F})$ arises locally as the Fefferman-type space of a Grassmann structure of type $(2, n)$; see Theorem 16 (ii).

Recall from Lemma 9 the explicit description of Weyl connections of an almost Grassmann structure and the one-to-one correspondence between Weyl structures and Weyl connections of an almost Grassmann structure. In particular, every Weyl connection is canonically associated to a collection of Weyl objects listed in (14). Now let $(M, E, F)$ be an almost Grassmann structure of type $(2, n)$, $n \geq 3$. Recall the convention of Penrose abstract index notation below (14). Given a Weyl connection $\nabla$ with canonically associated Weyl structure $\sigma$, denote by $R(\nabla^{TM})$, $R(\nabla^{E})$ and $R(\nabla^{F})$ the curvatures of $\nabla$ on $TM$, $E$ and $F$, respectively. Denote by $\tau$ the harmonic torsion and $\rho$ the harmonic curvature of $(M, E, F)$. Recall that $\tau$ agrees with the torsion of $\nabla$ on $TM$, see Lemma 9, and

$$\rho_{A'\,B'\,D'}^{A\;B\;C'} = R(\nabla^F)_{(A'\,B'\,D')}^{A\;B\;C'} - \frac{1}{n+2}\,R(\nabla^F)_{(A'\,B'|\,I'|}^{A\;B\;\;I'}\delta_{D')}^{C'} - \frac{2}{n+2}\,R(\nabla^F)^{[A\;\;B]I'}_{\;\;I'(A'\;\;B'}\delta_{D')}^{C'};$$

see, e.g., the end of Step (E) in [11, Subsection 4.1.3]. Moreover, denote by $Ric(\nabla)$ the Ricci curvature of $\nabla$ and denote by P and $(W, W')$ the Rho tensor and the Weyl tensor of $\sigma$, respectively. We prove the following relations.

**Lemma 18** *Notation as in the preceding paragraph. Then there holds*

$$\begin{cases} W_{A'\,B'\,D}^{A\;B\;C} &= R(\nabla^E)_{A'\,B'\,D}^{A\;B\;C} + \delta_D^B\,\mathrm{P}_{A'\,B'}^{A\;C} - \delta_D^A\,\mathrm{P}_{B'\,A'}^{B\;C} \\ W'^{A\;B\;C'}_{A'\,B'\,D'} &= R(\nabla^F)_{A'\,B'\,D'}^{A\;B\;C'} - \delta_{B'}^{C'}\,\mathrm{P}_{A'\,D'}^{A\;B} + \delta_{A'}^{C'}\,\mathrm{P}_{B'\,D'}^{B\;A} \end{cases} \tag{30}$$

*where*

$$\mathrm{P}_{A'\,B'}^{A\;B} = \frac{1}{n}\,Ric(\nabla)_{(A'\,B')}^{(A\;B)} + \frac{1}{n+4}\,Ric(\nabla)_{[A'\,B']}^{[A\;B]} + \frac{1}{n+2}\,Ric(\nabla)_{[A'\,B']}^{(A\;B)} + \frac{1}{n+2}\,Ric(\nabla)_{(A'\,B')}^{[A\;B]}. \tag{31}$$

***Proof*** Denote by $R(\nabla^{TM})$, $R(\nabla^{E})$ and $R(\nabla^{F})$ the curvatures of $\nabla$ on $TM$, $E$ and $F$, respectively. Since $TM \cong E^* \otimes F$, there holds

$$R(\nabla^{TM})_{A'\,B'\,C\;\;D'}^{A\;B\;C'D} = \delta_C^D\,R(\nabla^F)_{A'\,B'\,D'}^{A\;B\;C'} - \delta_{D'}^{C'}\,R(\nabla^E)_{A'\,B'\,C}^{A\;B\;D}$$

so that

$$\begin{aligned} Ric(\nabla)_{A'\,D'}^{A\;D} &= R(\nabla^{TM})_{A'\,I'\,I\;\;D'}^{A\;I\;I'D} \\ &= \delta_I^D\,R(\nabla^F)_{A'\,I'\,D'}^{A\;I\;I'} - \delta_{D'}^{I'}\,R(\nabla^E)_{A'\,I'\,I}^{A\;I\;D} \\ &= R(\nabla^F)_{A'\,I'\,D'}^{A\;DI'} - R(\nabla^E)_{A'\,D'\,I}^{A\;I\;D}. \end{aligned}$$

Denote by $(W, W')$ the Weyl tensor of $\sigma$ and by $\partial$ the Lie algebra homology differential; see (9). Recall from, e.g., [11, Theorem 5.2.3 (3)] that

$$(W, W') = (R(\nabla^E), R(\nabla^F)) + \partial\mathrm{P} \in \Omega^2(M, \mathfrak{s}(L(E, E) \oplus L(F, F))).$$

We write $\partial\mathrm{P} = (\Phi, \Psi)$ where $\Phi \in \Omega^2(M, L(E, E))$ and $\Psi \in \Omega^2(M, L(F, F))$. A computation shows that

$$(\Phi_{A'\,B'\,D}^{A\;B\;C}, \Psi_{A'\,B'\,D'}^{A\;B\;C'}) = (\delta_D^B\,\mathrm{P}_{A'\,B'}^{A\;C} - \delta_D^A\,\mathrm{P}_{B'\,A'}^{B\;C}, -\delta_{B'}^{C'}\,\mathrm{P}_{A'\,D'}^{A\;B} + \delta_{A'}^{C'}\,\mathrm{P}_{B'\,D'}^{B\;A}).$$

Hence there holds (30) and so

$$\begin{cases} W^{A\ I\ B}_{A'B'I} & = R(\nabla^E)^{A\ I\ B}_{A'B'I} + 2\,\mathrm{P}^{A\ B}_{A'B'} - \mathrm{P}^{A\ B}_{B'A'} \\ W'^{A\ B I'}_{A'I'B'} & = R(\nabla^F)^{A\ B I'}_{A'I'B'} - n\,\mathrm{P}^{A\ B}_{A'B'} + \mathrm{P}^{B\ A}_{A'B'}. \end{cases} \tag{32}$$

Recall from Lemma 8 that $W^{A\ I\ B}_{A'B'I} = W'^{A\ B I'}_{A'I'B'}$. Thus from (32) we obtain

$$Ric(\nabla)^{A\ B}_{A'B'} = (n+2)\,\mathrm{P}^{A\ B}_{A'B'} - \mathrm{P}^{A\ B}_{B'A'} - \mathrm{P}^{B\ A}_{A'B'}.$$

This equality is equivalent to (31) as each of them is equivalent to the following collection of relations.

$$Ric(\nabla)^{(A\ B)}_{(A'B')} = n\,\mathrm{P}^{(A\ B)}_{(A'B')} \qquad Ric(\nabla)^{[A\ B]}_{[A'B']} = (n+4)\,\mathrm{P}^{[A\ B]}_{[A'B']}$$
$$Ric(\nabla)^{(A\ B)}_{[A'B']} = (n+2)\,\mathrm{P}^{(A\ B)}_{[A'B']} \qquad Ric(\nabla)^{[A\ B]}_{(A'B')} = (n+2)\,\mathrm{P}^{[A\ B]}_{(A'B')};$$

cf. the end of [10, Subsection 4.3] □

We provide an equivalent description to Theorem 17.

**Corollary 19** *Let $(M, E, F)$ be an almost Grassmann structure of type $(2, n)$, $n \geq 3$. Denote by $\tau$ the harmonic torsion of $(M, E, F)$ and denote by $\rho$ the harmonic curvature of $(M, E, F)$. Given nowhere vanishing sections $\eta \in \Gamma(F)$ and $\varphi \in \Gamma(E^*)$, each of the following conditions is well-defined.*

- $(\nabla\eta)_o = 0$ *for one Weyl connection $\nabla$ and hence for all Weyl connections. Equivalently, $\eta$ is parallel for some Weyl connection.*
- $tr(\iota_\tau\tau)(\eta) = 0$ *where we view $tr(\iota_\tau\tau) \in \Gamma(F, T^*M \otimes E)$.*
- $(\nabla\varphi)_o = 0$ *for one Weyl connection $\nabla$ and hence for all Weyl connections. Equivalently, $\varphi$ is parallel for some Weyl connection.*
- $tr(\iota_\tau\tau)(\varphi) = 0$ *where we view $tr(\iota_\tau\tau) \in \Gamma(E^*, T^*M \otimes F^*)$.*
- $\frac{1}{2}\,\eta^{A'}\,\nabla^A_{A'}\varphi_A - \frac{1}{n}\,\varphi_A\,\nabla^A_{A'}\eta^{A'} = 1.$
- $\rho(\varphi \otimes \eta, TM) = 0.$
- $\tau(\varphi \otimes F, \varphi \otimes F) \subseteq \varphi \otimes F.$
- $\tau(\varphi \otimes \eta, TM) = 0$ *and, in addition, $(\varphi_A\,\eta^{B'}\,W^{A\ B\ C}_{A'B'D}, \varphi_A\,\eta^{B'}\,W'^{A\ B\ C'}_{A'B'D'}) = 0$ for the Weyl tensor $(W, W')$ of a Weyl structure and hence for the Weyl tensors of all Weyl structures.*

*Moreover, $(M, E, F)$ is locally the result of the Fefferman-type construction on a path geometry if and only if there exist nowhere vanishing sections $\eta \in \Gamma(F)$ and $\varphi \in \Gamma(E^*)$ such that all the above conditions hold.*

***Proof*** Denote by $\mathcal{T}$ the standard tractor of $(M, E, F)$. Then its dual $\mathcal{T}^*$ is the standard cotractor of $(M, E, F)$. Recall from [18, Proposition 5] that the trace-free part $(\nabla\eta)_o$ of $\nabla\eta \in \Gamma(E \otimes F^* \otimes F)$ is independent of the choice of Weyl connection $\nabla$ and the assignment $\eta \mapsto (\nabla\eta)_o$ is precisely the first Bernstein-Gelfand-Gelfand (BGG) operator of $\mathcal{T}$ and, similarly, the trace-free part $(\nabla\varphi)_o$ of $\nabla\varphi \in \Gamma(E \otimes F^* \otimes E^*)$ is independent of the choice of Weyl connection $\nabla$ and the assignment $\varphi \mapsto (\nabla\varphi)_o$ is precisely the first BGG operator of the standard cotractor bundle. Moreover, recall from [18, Proposition 5] that if $(\nabla\eta)_o = 0$, then $\eta$ is parallel for some Weyl connection and, similarly, $(\nabla\varphi)_o = 0$, then $\varphi$ is parallel for some Weyl connection. Thus we see that the first and the third conditions in the corollary are well-defined.

The symmetry (11) of $\tau$ implies that

$$\tau(\varphi \otimes \eta, TM) = 0 \quad \Longleftrightarrow \quad \tau(\varphi \otimes F, E^* \otimes \eta) = 0.$$

In this case, it is easy to see that $(W, W')(\varphi \otimes F, E^* \otimes \eta) = 0$ for the Weyl tensor $(W, W')$ of a Weyl structure implies that the same equality holds for the Weyl tensor $(W, W')$ of every Weyl structure. This shows that the last condition in the corollary is well-defined. That the other conditions make sense is obvious.

Recall from [18, Remark 6 and Theorem 8] the following. The natural projection $\mathcal{T} \twoheadrightarrow F$ induces a one-to-one correspondence

$$\{\mu \in \Gamma(\mathcal{T}) \text{ parallel}\} \longleftrightarrow \{\eta \in \Gamma(F) : (\nabla\eta)_o = 0 \text{ and } tr(\iota_\tau\tau)(\eta) = 0\}.$$

The inverse of this correspondence is the splitting operator $L : \Gamma(F) \to \Gamma(\mathcal{T})$. Similarly, the natural projection $\mathcal{T}^* \twoheadrightarrow E^*$ induces a one-to-one correspondence

$$\{\nu \in \Gamma(\mathcal{T}^*) \text{ parallel}\} \longleftrightarrow \{\varphi \in \Gamma(E^*) : (\nabla\varphi)_o = 0 \text{ and } tr(\iota_\tau\tau)(\varphi) = 0\}.$$

The inverse of this correspondence is the splitting operator $L : \Gamma(E^*) \to \Gamma(\mathcal{T}^*)$. Using the formulae of the two splitting operators in [18, Proposition 5], one computes that

$$(L\,\varphi)(L\,\eta) = \frac{1}{2}\,\eta^{A'}\,\nabla^{A}_{A'}\varphi_A - \frac{1}{n}\,\varphi_A\,\nabla^{A}_{A'}\eta^{A'}.$$

In particular, we see that the first five conditions in the corollary are equivalent to the conditions in Theorem 17 (i).

Denote by $\kappa$ the curvature of the normal Cartan geometry $(\mathcal{G} \to M, \omega)$ of type $(G, P)$ canonically associated to $(M, E, F)$. Recall from [5, Subsection 3.3 Proposition] that there holds $\kappa(\varphi \otimes \eta, TM) = 0$ if and only if both $\tau(\varphi \otimes \eta, TM) = 0$ and $\rho(\varphi \otimes \eta, TM) = 0$. In this case, there holds $\tau(\varphi \otimes F, E^* \otimes \eta) = 0$ so that, in particular, there holds $\kappa(\varphi \otimes \eta, TM) \subseteq \mathcal{G} \times_P \mathfrak{p}_+$ if and only if $(\varphi_A\,\eta^{B'}\,W^{A\ B\ C}_{A'B'D}, \varphi_A\,\eta^{B'}\,W'^{A\ B\ C'}_{A'B'D'}) = 0$ for the Weyl tensor $(W, W')$ of a Weyl structure and hence for the Weyl tensors of all Weyl structures. Finally, it is well-known that $\tau(\varphi \otimes F, \varphi \otimes F) \subseteq \varphi \otimes F$ if and only if $\varphi \otimes F$ is involutive; see, e.g., [11, Lemma 1.5.14]. We conclude that the last two conditions in the corollary are equivalent to the conditions in Theorem 17 (ii). This completes the proof of the corollary. □

In view of Lemma 18 and the paragraph preceding it, Corollary 19 is a characterization in terms of a Weyl connection and its associated tensors; it does not involve Cartan geometry. Moreover, each condition there is independent of the choice of Weyl connection.

### 5.1 The Case $n = 1$

The group homomorphism $i : G \to \tilde{G}$ given in (19) can be used to define a functor associating to every $(2n+1)$-dimensional generalized path geometry $(M, E, V)$ an almost Grassmann structure of type $(2, n+1)$. We have studied such constructions for all $n > 2$ as well as the restrictions of the construction over all geometries $(M, E, V)$ where $n = 2$ and $V$ is involutive, i.e., where $(M, E, V)$ is a path geometry. The cases where $n = 1$ and where $n = 2$ with $V$ non-involutive have been excluded from our previous discussion due to their different types of fundamental invariants; see the beginning of the article. We discuss the case where $n = 1$ below.

First, we recall the following background from, e.g., [11, Subsections 4.1.3, 4.1.4 and 4.2.3]. In dimension $2 \cdot 1 + 1 = 3$, the concepts of generalized path geometry and of Lagrangian contact structure coincide. Such a geometry $(M, E, V)$ is automatically torsion-free. Its fundamental invariants are given by two Cotten-York-type tensors in

$$\Gamma(E^* \otimes (TM/(E \oplus V))^* \otimes E^*) \quad \text{and} \quad \Gamma(V^* \otimes (TM/(E \oplus V))^* \otimes V^*), \tag{33}$$

respectively. On the other hand, a nonzero volume form on $\mathbb{R}^4$ gives rise to a natural identification $\wedge^2 \mathbb{R}^4 \cong \mathbb{R}^{(3,3)}$, which in turn induces an isomorphism

$$\tilde{G} = SL(4, \mathbb{R}) \cong Spin(3, 3).$$

It restricts to an isomorphism between the structure groups

$$\tilde{P}_0 \cong CSpin(2, 2).$$

Here, $\tilde{P}_0$ denotes the Levi subgroup of $\tilde{P}$, that is, the block diagonal subgroup of $SL(4, \mathbb{R})$ where the block size is $(2, 2) \times (2, 2)$. In particular, an almost Grassmann structure $(\tilde{M}, \tilde{E}, \tilde{F})$ of type $(2, 2)$ is the same as a conformal spin structure of signature $(2, 2)$. Moreover, $(\tilde{M}, \tilde{E}, \tilde{F})$ is automatically torsion-free. Its fundamental invariants are a pair of harmonic curvatures in

$$\Gamma(\wedge^2 \tilde{E} \otimes \mathrm{Sym}^2 \tilde{F}^* \otimes \mathfrak{sl}(\tilde{F})) \quad \text{and} \quad \Gamma(\mathrm{Sym}^2 \tilde{E} \otimes \wedge^2 \tilde{F}^* \otimes \mathfrak{sl}(\tilde{E})),$$

respectively.

Recall in addition that the spin representation of $Spin(3, 3)$ corresponds to the natural representation of $SL(4, \mathbb{R})$ on $\mathbb{R}^4 \oplus \mathbb{R}^{4*}$ and the Clifford multiplication decomposes into the maps

$$\wedge^2 \mathbb{R}^4 \times \mathbb{R}^4 \to \mathbb{R}^{4*} \quad \text{and} \quad \wedge^2 \mathbb{R}^4 \times \mathbb{R}^{4*} \to \mathbb{R}^4$$

canonically induced by the volume form on $\mathbb{R}^4$ chosen before. In particular, every spinor $0 \neq s \in \mathbb{R}^4 \oplus \mathbb{R}^{4*}$ is pure, that is, it is annihilated by a maximal isotropic subspace of $\wedge^2 \mathbb{R}^4$. Moreover, the almost Grassmann standard tractor bundle $\tilde{\mathcal{T}}$ is the positive conformal spin tractor bundle and its dual $\tilde{\mathcal{T}}^*$ is the negative conformal spin tractor bundle.

The authors of [21] introduced a Fefferman-type construction based on a group homomorphism

$$SL(n+2,\mathbb{R}) \to Spin(n+2,n+2)$$

that associates to every $(2n+1)$-dimensional Lagrangian contact structure a conformal spin structure of signature $(n+1,n+1)$; see [21, Subsection 3.2]. In the case where $n=1$, the group homomorphism introduced in [21, Subsection 3.2] coincides with our homomorphism

$$i : SL(3,\mathbb{R}) \to SL(4,\mathbb{R}) \cong Spin(3,3);$$

see (19). In particular, the two Fefferman-type constructions are the same in this dimension up to minor differences on the choice of parabolic subgroups $Q \subseteq G$ and $\tilde{P} \subseteq \tilde{G}$.

We obtain the following results, which agree with [21, Proposition 3.7 and Proposition 3.10] under the setup of [21, Subsection 3.2] with $n=2$ in their work. As mentioned above, every tractor spinor is pure in this specific dimension.

**Theorem 20** *The Fefferman-type construction from 3-dimensional generalized path geometries to almost Grassmann structures of type (2, 2) determined by (19) is normal.*

*Moreover, an almost Grassmann structure $(\tilde{M}, \tilde{E}, \tilde{F})$ of type (2, 2) is locally the result of the Fefferman-type construction on a 3-dimensional generalized path geometry if and only if the conditions in Theorem 17 (i) are satisfied. In terms of conformal spin structure, these conditions reads as follows:*

*There are pure parallel tractor spinors $\mu \in \Gamma(\tilde{\mathcal{T}})$ and $\nu \in \Gamma(\tilde{\mathcal{T}}^*)$ with $\mu(\nu) = 1$. Moreover, the projection $\eta \in \Gamma(\tilde{F})$ of $\mu$ and $\varphi \in \Gamma(\tilde{E}^*)$ of $\nu$ to the underlying twistor spinors are both nowhere vanishing.*

***Proof*** Since the two Cotten-York-type invariants (33) above are still sections of the $Q$-module $\mathbb{F}$ defined in the proof of Lemma 14 (i), the statement of Lemma 14 (i) holds for $n=1$ as well. Since the proof of Proposition 15 is a computation based on Lemma 14 (i), it applies to the case $n=1$. Hence the statement of Proposition 15 holds true for $n=1$ too. This proves the first assertion of the theorem.

The proof of Theorem 17 can be directly applied to the case $n=1$. Since $(\tilde{M}, \tilde{E}, \tilde{F})$ is automatically torsion-free, Theorem 17 gives the second assertion of the theorem. □

The Fefferman-type construction induced by the related homomorphism

$$SL(n+2,\mathbb{R}) \to SO(n+2,n+2)$$

associating to every Lagrangian contact structure a conformal structure has been studied in [24]. Note that

$$Spin(n+2,n+2) \to SO(n+2,n+2)$$

is a two-fold covering. In particular, we obtain the following from the last relation in the proof of [24, Theorem 3.3]. Given a 3-dimensional generalized path geometry $(M, E, V)$ with canonically associated normal regular parabolic geometry $(\mathcal{G} \to M, \omega)$, there exists local coordinates $(x, u, p) = U \subseteq M$ and a local section $\phi \in \Gamma_U(\mathcal{G})$ such that

$$\phi^*\omega = \begin{pmatrix} * & * & * \\ \theta & \alpha & * \\ \sigma & \pi & * \end{pmatrix} \in \Omega^1(U, \mathfrak{sl}(3, \mathbb{R})).$$

Here,

$$\sigma = du - p\,dx, \quad \theta = dx, \quad \pi = dp - f\,dx, \quad \alpha = \frac{1}{6}\frac{\partial^2 f}{\partial p^2}\sigma + \frac{2}{3}\frac{\partial f}{\partial p}\theta$$

and $f = f(x, u, p)$ is the defining function of the geometry with respect to the given local coordinates; see [15, (2.1)]. On the other hand, the projection of $\phi \in \Gamma_U(\mathcal{G})$ to $\underline{\phi} \in \Gamma_U(\tilde{M})$ defines a local trivialization of the principal bundle $\tilde{M} = GL_+(E^* \otimes V)$. Note that the isomorphism between $Q/H$ and the structure group $\mathbb{R}_+$ of $GL_+(E^* \otimes V)$ is defined by the surjective homomorphism

$$Q \to \mathbb{R}_+, \quad \begin{pmatrix} b & * & * \\ 0 & a & * \\ 0 & 0 & (ab)^{-1} \end{pmatrix} \mapsto a^{-3}$$

with kernel $H$. Let $\tilde{s} : \mathbb{R}_+ \to \mathbb{R}$ be the inverse of the Lie group isomorphism

$$\mathbb{R} \xrightarrow{\cong} \mathbb{R}_+, \quad X \mapsto e^{-3X}.$$

Thus $(x, u, p, \tilde{s})$ are local coordinates on the preimage $\tilde{U} \subseteq \tilde{M}$ of $U \subseteq M$. Note that $\tilde{s}$ is related to $s$ defined in Subsection 3.3 by $\tilde{s} = -2s$. In particular, in the spirit of (3.5), consider the local section $\psi : \tilde{U} \to \mathcal{G}$ defined by

$$\psi(\underline{\phi}(y)\, e^{-3X}) = \phi(y) \begin{pmatrix} e^{-\frac{1}{2}X} & 0 & 0 \\ 0 & e^X & 0 \\ 0 & 0 & e^{-\frac{1}{2}X} \end{pmatrix}$$

for all $y \in U$ and $X \in \mathbb{R}$. There holds

$$\psi^*\omega(\frac{\partial}{\partial \tilde{s}}) = \begin{pmatrix} -\frac{1}{2} & 0 & 0 \\ 0 & 1 & 0 \\ 0 & 0 & -\frac{1}{2} \end{pmatrix}.$$

We also denote the pullbacks of $\sigma, \theta, \pi, \alpha$ along $\tilde{U} \to U$ by the same letters, which are by definition invariant one-forms on $\tilde{U}$. Moreover, observe that $d\tilde{s} \in \Omega^1(\tilde{U})$ is also invariant. In particular, on $\underline{\phi}(U)$, there holds

$$\psi^*\omega|_{\underline{\phi}(U)} = \begin{pmatrix} * & * & * \\ \theta & \alpha + d\tilde{s} & * \\ \sigma & \pi & * \end{pmatrix}.$$

By equivariancy of $\omega$, there holds

$$\psi^*\omega = \begin{pmatrix} * & * & * \\ e^{-\frac{3\tilde{s}}{2}}\theta & \alpha + d\tilde{s} & * \\ \sigma & e^{\frac{3\tilde{s}}{2}}\pi & * \end{pmatrix} \in \Omega^1(\tilde{U}, \mathfrak{sl}(3, \mathbb{R})).$$

Hence

$$\begin{pmatrix} e^{-\frac{3\tilde{s}}{2}}\theta & \alpha + d\tilde{s} \\ \sigma & e^{\frac{3\tilde{s}}{2}}\pi \end{pmatrix}$$

or, equivalently,

$$\begin{pmatrix} \theta & \alpha + d\tilde{s} \\ \sigma & \pi \end{pmatrix} = \begin{pmatrix} e^{\tilde{s}} & 0 \\ 0 & e^{-\frac{\tilde{s}}{2}} \end{pmatrix} \begin{pmatrix} e^{-\frac{3\tilde{s}}{2}}\theta & \alpha + d\tilde{s} \\ \sigma & e^{\frac{3\tilde{s}}{2}}\pi \end{pmatrix} \begin{pmatrix} e^{\frac{\tilde{s}}{2}} & 0 \\ 0 & e^{-\tilde{s}} \end{pmatrix}$$

is a local coframe on $\tilde{M}$ adapted to the almost Grassmann structure $(\tilde{M}, \tilde{E}, \tilde{F})$. Moreover, the metric in the conformal class of the geometry is given by

$$\begin{aligned} \tilde{g} &= \theta \odot \pi - \sigma \odot (\alpha + d\,\tilde{s}) \\ &= \frac{1}{2}(\theta \otimes \pi + \pi \otimes \theta - \sigma \otimes (\alpha + d\,\tilde{s}) - (\alpha + d\,\tilde{s}) \otimes \sigma). \end{aligned}$$

This agrees with the formula of the metric provided in Theorem 3.3 and (6.4). Note that our $d\,\tilde{s}$ corresponds to $-2\,d\,s$ there.

## 6 The Fefferman-Type Construction on Almost Grassmann Structures

We follow Convention 6 for the notation of Lie groups, Lie algebras and gradings. Recall that every almost Grassmann structure

$$(M, E, F)$$

of type $(2, n)$ is canonically associated to a normal parabolic geometry

$$(\mathcal{G} \to M, \omega)$$

of type $(G, P)$. Replacing $Q$ by $P$ in the argument in the first paragraph of Sect. 5, we see that the Fefferman-type functor defined by $i : G \to \tilde{G}$ in (19) associates to

$(\mathcal{G} \to M, \omega)$ a parabolic geometry

$$(\tilde{\mathcal{G}} \to \tilde{M}, \tilde{\omega})$$

of type $(\tilde{G}, \tilde{P})$ with underlying almost Grassmann structure

$$(\tilde{M}, \tilde{E}, \tilde{F})$$

of type $(2, n + 1)$. This defines a functorial construction associating to every almost Grassmann structure of type $(2, n)$ an almost Grassmann structure of type $(2, n + 1)$. We denote by $\kappa$ and $\tilde{\kappa}$ the curvatures of $\omega$ and $\tilde{\omega}$, respectively. Denote by

$$j : \mathcal{G} \to \tilde{\mathcal{G}}$$

the natural map defining $\tilde{\mathcal{G}}$ and denote by

$$\pi : \tilde{\mathfrak{g}}/\tilde{\mathfrak{p}} \twoheadrightarrow \mathfrak{g}/\mathfrak{p}$$

the projection canonically induced by the Lie algebra map $i'$ of $i : G \to \tilde{G}$; see Sect. 3.

We provide a geometric description of $\tilde{M}$, $\tilde{E}$ and $\tilde{F}$.

**Theorem 21** *Notation as above. The following canonical identifications hold.*

- *The fiber bundle $p : \tilde{M} \to M$ is isomorphic to $E^* \setminus \{0\}$, i.e., the complement of the zero section in $E^*$.*
- *$\tilde{E}$ is the pullback bundle $p^*E$. Additionally, $\tilde{E}^*$ is the vertical bundle of $p : \tilde{M} \to M$.*
- *$\tilde{F} = (p^*\mathcal{T})/\ell$. Here, $p^*\mathcal{T}$ is the pullback of the standard tractor bundle $\mathcal{T}$ of $(M, E, F)$ and $\ell$ is the tautological subbundle defined by*

$$\ell(\varphi_x) = \ker \varphi_x \subseteq E_x \subseteq \mathcal{T}_x$$

*over every $\varphi_x \in E_x^* \setminus \{0\} \subseteq \tilde{M}$ for each $x \in M$.*

***Proof*** Notation as at the beginning of the section. Consider the natural representation of $P$ on $\mathbb{R}^{2*}$. The isotropy subgroup fixing the base point $(0, 1) \in \mathbb{R}^{2*}$ is precisely $H$ and the $P$-orbit of the base point is $\mathbb{R}^{2*} \setminus \{0\}$. Since $E^* = \mathcal{G} \times_P \mathbb{R}^{2*}$,

$$\tilde{M} = \mathcal{G}/H = \mathcal{G} \times_P P/H = E^* \setminus \{0\}$$

as fiber bundles over $M$. This proves the first assertion of the theorem. Note that, in contrast to Theorem 11, there is no natural principal bundle structure on $p : \tilde{M} \to M$ because $H$ is not a normal subgroup of $P$.

Using the same argument from the proof of Theorem 11, we obtain the relations

$$\tilde{E} = \mathcal{G} \times_H \mathbb{R}^2 \quad \text{and} \quad \tilde{F} = \mathcal{G} \times_H \mathbb{R}^{n+2}/\mathbb{R}$$

where

$$\mathbb{R} = (\mathbb{R}, 0)^t \subseteq \mathbb{R}^2 \subseteq \mathbb{R}^{n+2}$$

is the filtration determined by the natural representation of $H$ on $\mathbb{R}^{n+2}$. In particular, $\tilde{E} = p^*E$. Moreover, the isomorphism $\mathbb{R}^{2*} \cong \mathfrak{p}/\mathfrak{h}$ of $H$-modules implies that $\tilde{E}^*$ is canonically identified with the vertical distribution $\mathcal{G} \times_H \mathfrak{p}/\mathfrak{h} \subseteq T\tilde{M}$. This completes the proof of the second assertion of the theorem.

Finally, there holds $\tilde{F} = (p^*\mathcal{T})/\ell$ because $\mathcal{G} \times_H \mathbb{R}^{n+2} = p^*\mathcal{T}$ and $\mathcal{G} \times_H (\mathbb{R}, 0)^t = \ell$. This proves the final assertion. □

***Remark 22*** Observe that the sections

$$\varphi = \mathcal{G} \times_H (0, 1) \in \Gamma(\tilde{E}^*) \quad \text{and} \quad \eta = \mathcal{G} \times_H (-1, 0^n)^t \in \Gamma(\tilde{F})$$

we specified in Remark 12 are also well-defined in the above situation. Here, $\varphi$ is the tautological section of $\tilde{E}^*$ so that its rank-1-kernel is precisely the line bundle

$$\ell \subseteq \tilde{E}$$

defined in Theorem 21. Moreover, the vertical distribution of $p : \tilde{M} \to M$ is

$$\tilde{V} = \tilde{E}^* \otimes \eta \subseteq T\tilde{M}.$$

**Convention 23** *We define a decomposition*

$$\tilde{\mathfrak{g}} = \mathfrak{sl}(n+3, \mathbb{R}) = \left(\begin{array}{c|c|c} \tilde{\mathfrak{n}}_0^E & \tilde{\mathfrak{n}}_1^E & \tilde{\mathfrak{n}}_2 \\ \hline \tilde{\mathfrak{n}}_{-1}^E & \tilde{\mathfrak{n}}_0^\eta & \tilde{\mathfrak{n}}_1^F \\ \hline \tilde{\mathfrak{n}}_{-2} & \tilde{\mathfrak{n}}_{-1}^F & \tilde{\mathfrak{n}}_0^F \end{array}\right)$$

*in the block size* $(2, 1, n) \times (2, 1, n)$. *Observe that it is a grading decomposition of some parabolic subgroup of* $\tilde{G}$ *that contains* $H$. *In particular, its induced filtration is* $H$*-invariant and each block component is invariant under the natural representation of*

$$H_0 = P_0 \cap H = \tilde{P}_0 \cap i(H). \tag{34}$$

*Observe that*

$$\tilde{\mathfrak{n}}^{1,F} = \tilde{\mathfrak{n}}_1^F \oplus \tilde{\mathfrak{n}}_2$$

*is an* $H$*-submodule of* $\tilde{\mathfrak{g}}$.

We say that the Fefferman-type construction on an almost Grassmann structure $(M, E, F)$ of type $(2, n)$ is normal if the Cartan connection $\tilde{\omega}$ given at the beginning of the section is normal. We prove that this notion of being normal is equivalent to a certain property of the torsion $\tau$ of $(M, E, F)$. Recall from Lemma 8 the contraction $tr(\iota_\tau \tau) \in \Gamma(\mathrm{Sym}^2 T^*M)$ of $\tau \otimes \tau$.

**Proposition 24** *The Fefferman-type construction on $(M, E, F)$ is normal if and only if there holds $tr(\iota_\tau \tau) = 0$.*

***Proof*** Notation as at the beginning of the section. For every $\tilde{x} \in \tilde{M}$, we may choose $u \in \mathcal{G}$ such that $j(u) \in \tilde{\mathcal{G}}$ projects to $\tilde{x}$. Denote by $x \in M$ the projection of $\tilde{x}$ along $p : \tilde{M} \to M$. Choose a Weyl structure $\sigma : \mathcal{G}_0 = \mathcal{G}/P_+ \to \mathcal{G}$ whose image contains $u$ so that, in particular, we consider the identification between $\mathfrak{g}$ and the fiber of $\mathcal{G} \times_P \mathfrak{g}$ over $x \in M$ defined by $u$. Thus

$$\kappa(u) = \left(\begin{array}{c|c} W(x)_{A'\,B'\,D}^{A\;\,B\;\,C} & Y(x)_{A'\,B'\,D'}^{A\;\,B\;\,C} \\ \hline \tau(x)_{A'\,B'\,D}^{A\;\,B\;\,C'} & W'(x)_{A'\,B'\,D'}^{A\;\,B\;\,C'} \end{array}\right). \tag{35}$$

Here, $(W, W')$ is the Weyl tensor of $\sigma$ and $Y$ is the Cotten-York tensor of $\sigma$; see (14). These Weyl objects evaluated at $x$ are independent of the choice of the Weyl structure mentioned above. We view

$$\varphi = \mathcal{G} \times_H (0, 1) \in \Gamma(\tilde{E}^*)$$

defined in Remark 22 as a section of $p^*E^* = \tilde{E}^*$; see Theorem 21. Using the identification between $\tilde{\mathfrak{g}}$ and the fiber of $\tilde{\mathcal{G}} \times_{\tilde{P}} \tilde{\mathfrak{g}}$ over $\tilde{x} \in \tilde{M}$ defined by $j(u) \in \tilde{\mathcal{G}}$, there holds

$$\tilde{\kappa}(j(u)) = \left(\begin{array}{c|c|c} W(x)_{A'\,B'\,D}^{A\;\,B\;\,C} & 0 & Y(x)_{A'\,B'\,D'}^{A\;\,B\;\,C} \\ \hline \varphi(\tilde{x})_I\, W(x)_{A'\,B'\,D}^{A\;\,B\;\,I} & 0 & \varphi(\tilde{x})_I\, Y(x)_{A'\,B'\,D'}^{A\;\,B\;\,I} \\ \hline \tau(x)_{A'\,B'\,D}^{A\;\,B\;\,C'} & 0 & W'(x)_{A'\,B'\,D'}^{A\;\,B\;\,C'} \end{array}\right) \circ \wedge^2 \pi.$$

The decomposition of $\tilde{\mathfrak{g}}$ is as given in Convention 23.

Denote by $\tilde{\partial}^*$ the Kostant codifferential of $(\tilde{\mathfrak{g}}, \tilde{\mathfrak{p}})$. We show that

$$\frac{1}{2}\,(\tilde{\partial}^*\tilde{\kappa})(j(u)) = \left(\begin{array}{c|c|c} 0 & 0 & 0 \\ \hline 0 & 0 & -\varphi_J(\tilde{x})\, W(x)_{A'C'I}^{A\;\,I\;\,J} \\ \hline 0 & 0 & 0 \end{array}\right) \circ \pi. \tag{36}$$

Indeed, the dual of $\pi : \tilde{\mathfrak{g}}/\tilde{\mathfrak{p}} \twoheadrightarrow \mathfrak{g}/\mathfrak{p}$ is

$$\pi^* : \mathfrak{p}_+ \to \tilde{\mathfrak{p}}_+, \quad \pi^*(Z)_{ij} = Z_{i,j+1}, \quad Z \in \mathfrak{p}_+.$$

Note that it coincides with the restriction of $\pi^*$ defined in the proof of Proposition 15. We continue to denote by $\alpha : \mathfrak{g} \to \tilde{\mathfrak{g}}$ the embedding of $\mathfrak{g}$ into the Lie subalgebra of $\tilde{\mathfrak{g}}$ with zeros in the third row and the third column. Then

$$\left[\pi^*\left(\begin{array}{c|c} 0 & Z_{12} \\ \hline 0 & 0 \end{array}\right), i'\left(\begin{array}{c|c} \Phi_{11} & \Phi_{12} \\ \hline \Phi_{21} & \Phi_{22} \end{array}\right)\right] = \alpha([Z, \Phi]) - \left(\begin{array}{c|c|c} 0 & 0 & 0 \\ \hline 0 & 0 & (0,1)\,\Phi_{11} Z_{12} \\ \hline 0 & 0 & 0 \end{array}\right)$$

for all $Z = Z_{12} \in \mathfrak{p}_+$ and $\Phi = (\Phi_{ij})_{1\le i,j\le 2} \in \mathfrak{g}$. The block size on the left-hand side is $(2, n) \times (2, n)$ and the block size on the right-hand side is $(2, 1, n) \times (2, 1, n)$. Following steps analogous to those succeeding (25), one obtains the expression above.

It follows from Lemma 8 that $\varphi_J(\tilde{x})\, W(x)_{A'C'I}^{A\ I\ J} = 0$ if and only if $\varphi_J(\tilde{x})\, tr(\iota_\tau \tau)(x)_{A'C'}^{A\ J} = 0$, that is, there holds $(\tilde{\partial}^*\tilde{\kappa})(\tilde{x}) = 0$ if and only if $\varphi_J(\tilde{x})\, tr(\iota_\tau \tau)(x)_{A'C'}^{A\ J} = 0$. Since the $P$-orbit of $(0, 1)$ is $\mathbb{R}^{2*} \setminus \{0\}$, given any nonzero element $v \in E_x^*$, there is a unique point $\tilde{x} \in \tilde{M}$ over $x$ such that $\varphi(\tilde{x}) = v$ and hence there holds $\tilde{\partial}^*\tilde{\kappa} = 0$ if and only if $tr(\iota_\tau \tau) = 0$. This proves the statement of the proposition. □

Recall from, e.g., [11, Proposition 3.1.13] that we may normalize $(\tilde{\mathcal{G}} \to \tilde{M}, \tilde{\omega})$ and thus obtain the normal parabolic geometry $(\tilde{\mathcal{G}} \to \tilde{M}, \tilde{\omega}^{nor})$ canonically associated to $(\tilde{M}, \tilde{E}, \tilde{F})$. Denote the curvature of $\tilde{\omega}^{nor}$ by $\tilde{\kappa}^{nor}$. Let the group $H_0$ and the grading decomposition of $\tilde{\mathfrak{g}}$ be as in Convention 23. Moreover, let $\tilde{V} \subseteq T\tilde{M}$ be the vertical distribution of the projection $\tilde{M} \to M$ as in Remark 22.

**Lemma 25** *Notation as above. Then the following properties hold.*

- $j^*(\tilde{\omega}^{nor} - \tilde{\omega}) \in \Omega^1(\mathcal{G}, \tilde{\mathfrak{n}}^{1,F})^H$.
- $(\tilde{\omega}^{nor} - \tilde{\omega})(\tilde{V}) = 0$.
- $j^*(\tilde{\kappa}^{nor} - \tilde{\kappa}) \in \Omega^2(\mathcal{G}, [\tilde{\mathfrak{g}}, \tilde{\mathfrak{n}}^{1,F}])^H$.
- $(j^*\tilde{\kappa}^{nor})(\tilde{V}, T\tilde{M}) \subseteq \mathcal{G} \times_H \tilde{\mathfrak{n}}^{1,F}$.
- $\tilde{\kappa}^{nor}(\tilde{V}, \tilde{V}) = 0$.

***Proof*** Consider the $H$-modules

$$\begin{aligned}
\mathbb{E} &= \tilde{\mathfrak{n}}_2 \otimes \tilde{\mathfrak{n}}^{1,F} \\
&= \{\phi \in \tilde{\mathfrak{p}}_+ \otimes \tilde{\mathfrak{n}}^{1,F} : \phi(i'(\mathfrak{p})) = 0\} \quad \text{and} \\
\mathbb{F} &= \tilde{\mathfrak{n}}_1^E \wedge \tilde{\mathfrak{n}}_2 \otimes \tilde{\mathfrak{n}}^{1,F} \oplus \wedge^2 \tilde{\mathfrak{n}}_2 \otimes [\tilde{\mathfrak{g}}, \tilde{\mathfrak{n}}^{1,F}] \\
&= \{\psi \in \wedge^2 \tilde{\mathfrak{p}}_+ \otimes [\tilde{\mathfrak{g}}, \tilde{\mathfrak{n}}^{1,F}] : \phi(i'(\mathfrak{p}), i'(\mathfrak{p})) = 0 \text{ and } \phi(i'(\mathfrak{p}), i'(\mathfrak{g})) \subseteq \tilde{\mathfrak{n}}^{1,F}\}.
\end{aligned}$$

Note that the indicated decomposition of $\mathbb{F}$ is only $H_0$-invariant. Let $\tilde{\partial}$ be the Lie algebra differential and $\tilde{\partial}^*$ the Kostant codifferential of $(\tilde{\mathfrak{g}}, \tilde{\mathfrak{p}})$. A computation shows that

$$\tilde{\partial}\,\mathbb{E} \subseteq \mathbb{F} \quad \text{and} \quad \tilde{\partial}^*\,\mathbb{F} \subseteq \mathbb{E} \tag{37}$$

so that, in particular, it follows from the Hodge decomposition (10) that each of $\tilde{\partial}$ and $\tilde{\partial}^*$ induces a linear isomorphism

$$\operatorname{im}\tilde{\partial}^* \cap \mathbb{E} \underset{\tilde{\partial}^*}{\overset{\tilde{\partial}}{\rightleftharpoons}} \operatorname{im}\tilde{\partial} \cap \mathbb{F}. \tag{38}$$

in opposite directions. Here, $\tilde{\partial}$ is $H_0$-equivariant and $\tilde{\partial}^*$ is $H$-equivariant.

On the other hand, the expression (36) in the proof of Proposition 24 implies that

$$j^*(\tilde{\partial}^*\tilde{\kappa}) \in \Gamma(\mathcal{G} \times_H \mathbb{E}). \tag{39}$$

Since $\mathbb{E}^{(1)} = \mathbb{E} \subseteq \tilde{\mathfrak{p}}_1 \otimes \tilde{\mathfrak{p}}^0$, we define

$$\mathbb{E}^{(2)} = \mathbb{E} \cap (\tilde{\mathfrak{p}}_1 \otimes \tilde{\mathfrak{p}}_1) = \tilde{\mathfrak{n}}_2 \otimes \tilde{\mathfrak{n}}_2,$$

which is an $H$-submodule of $\mathbb{E}$. We follow the normalizing procedure in the proof of [11, Proposition 3.1.13] in our specific setting. Indeed, in view of (38), for every local trivialization $\tilde{U} \times H \subseteq \mathcal{G}$ with $\tilde{U} \subseteq \tilde{M}$ open, we may choose an $H$-equivariant map $\phi : \tilde{U} \times H \to \mathbb{E}$ such that the image of the $H$-equivariant map

$$\tilde{\partial}^*\tilde{\partial}\phi + j^*(\tilde{\partial}^*\tilde{\kappa}) : \tilde{U} \times H \to \mathbb{E}$$

is contained in $\mathbb{E}^{(2)}$. Using partition of unity, we obtain $\phi \in \Gamma(\mathcal{G} \times_H \mathbb{E})$ such that

$$\tilde{\partial}^*\tilde{\partial}\phi + j^*(\tilde{\partial}^*\tilde{\kappa}) \in \Gamma(\mathcal{G} \times_H \mathbb{E}^{(2)}).$$

$\phi$ can be extended uniquely to $\Gamma(\tilde{\mathcal{G}} \times_{\tilde{P}} (\tilde{\mathfrak{p}}_1 \otimes \tilde{\mathfrak{p}}^0))$. Now define a Cartan connection

$$\tilde{\omega}' = \tilde{\omega} + \phi \circ \tilde{\omega}$$

with curvature $\tilde{\kappa}'$. A computation using the formula of Cartan curvature in Definition 1 shows that

$$j^*(\tilde{\kappa}' - \tilde{\kappa}) \in \Gamma(\mathcal{G} \times_H \mathbb{F}).$$

Moreover, there holds

$$j^*(\tilde{\partial}^*\tilde{\kappa}') \in \Gamma(\mathcal{G} \times_H \mathbb{E}^{(1)}).$$

Repeating the above procedure once more, we obtain the normalized connection $\tilde{\omega}^{nor}$ which satisfies the following properties.

$$j^*(\tilde{\omega}^{nor} - \tilde{\omega}) \in \Gamma(\mathcal{G} \times_H \mathbb{E}) \quad \text{and} \quad j^*(\tilde{\kappa}^{nor} - \tilde{\kappa}) \in \Gamma(\mathcal{G} \times_H \mathbb{F}). \tag{40}$$

We obtain all assertions of the lemma by putting (40), the properties

$$j^*\tilde{\omega} \in \Omega^1(\mathcal{G}, i'(\mathfrak{g})) \quad \text{and} \quad \tilde{\kappa}(\tilde{V}, T\tilde{M}) = 0$$

of the Fefferman-type construction together with the relations

$$\tilde{\mathfrak{n}}^{1,F} \subseteq [\tilde{\mathfrak{g}}, \tilde{\mathfrak{n}}^{1,F}] \quad \text{and} \quad [i'(\mathfrak{p}), \tilde{\mathfrak{n}}^{1,F}] \subseteq \tilde{\mathfrak{n}}^{1,F}.$$

This completes the proof of the lemma. □

Consider the standard tractor bundle $\tilde{\mathcal{T}} = \tilde{\mathcal{G}} \times_{\tilde{P}} \mathbb{R}^{n+3}$. Let $\eta \in \Gamma(\tilde{F})$, $\varphi \in \Gamma(\tilde{E}^*)$ and $\tilde{V} \subseteq T\tilde{M}$ be as in Remark 22. We obtain a refined filtration

$$\tilde{E} \subseteq \operatorname{span}(\eta) + \tilde{E} \subseteq \tilde{\mathcal{T}}$$

of ranks $2 \leq 3 \leq n+3$. Its dual filtration is

$$\operatorname{Ann}(\eta) \subseteq \tilde{F}^* \subseteq \tilde{\mathcal{T}}^*$$

with ranks $n \leq n+1 \leq n+3$. Here, $\operatorname{Ann}(\eta)$ denotes the annihilators of $\eta$ in $\tilde{F}^*$. Both the standard tractor bundle $\tilde{\mathcal{T}}$ and the standard cotractor bundle $\tilde{\mathcal{T}}^*$ are equipped with a tractor connection induced by $\tilde{\omega}^{nor}$. We denote them by $\nabla^{\tilde{\omega}^{nor}}$.

**Corollary 26** *There is a unique standard tractor $\mu \in \Gamma(\tilde{\mathcal{T}})$ with projection $\eta \in \Gamma(\tilde{F})$ and a unique standard cotractor $\nu \in \Gamma(\tilde{\mathcal{T}}^*)$ with projection $\varphi \in \Gamma(\tilde{E}^*)$ such that $\nu(\mu) = 1$, $\nabla^{\tilde{\omega}^{nor}} \mu = 0$ and that*

$$\nabla^{\tilde{\omega}^{nor}} \nu \in \{\Phi \in \Omega^1(\tilde{M}, \operatorname{Ann}(\eta)) \subseteq \Omega^1(\tilde{M}, \tilde{F}^*) : \Phi(\tilde{V}) = 0\}.$$

*Moreover, there holds*

$$\tilde{\kappa}^{nor}(\tilde{V}, T\tilde{M}) \subseteq \tilde{\mathcal{N}}^{1,F}$$

*where*

$$\tilde{\mathcal{N}}^{1,F} = \{\Psi \in \Gamma(\mathfrak{sl}(\tilde{\mathcal{T}})) : \Psi(\operatorname{span}(\eta) + \tilde{E}) = 0 \text{ and } \Psi(\tilde{\mathcal{T}}) \subseteq \operatorname{span}(\eta) + \tilde{E}\}.$$

***Proof*** The uniqueness of $\mu$ and $\nu$ follows from general theory; see Remark 28 below. We now prove the existence of $\mu$ and $\nu$. To see this, observe that the argument in the proof of Theorem 16 (iii) applies here. Precisely, the sections

$$\mu = \mathcal{G} \times_H (0, 0, -1, 0^n)^t \in \Gamma(\tilde{\mathcal{T}}) \quad \text{and} \quad \nu = \mathcal{G} \times_H (0, 1, -1, 0^n) \in \Gamma(\tilde{\mathcal{T}}^*)$$

are well-defined, project to $\eta$ and $\varphi$, respectively, and such that $\nu(\mu) = 1$. Moreover, denote by $\nabla^{\tilde{\omega}}$ the tractor connections on $\tilde{\mathcal{T}}$ and on $\tilde{\mathcal{T}}^*$ induced by the Cartan connection $\tilde{\omega}$ as given at the beginning of the section. We have proven that

$$\nabla^{\tilde{\omega}} \mu = 0 \quad \text{and} \quad \nabla^{\tilde{\omega}} \nu = 0.$$

On the other hand, observe that the relation (6) in the proof of Theorem 4 applies. Together with the properties of $\tilde{\omega}^{nor} - \tilde{\omega}$ in Lemma 25, we compute that

$$\begin{aligned}
&(\nabla^{\tilde{\omega}^{nor}}_{\xi}\mu - \nabla^{\tilde{\omega}}_{\xi}\mu)(u) \subseteq \tilde{\mathfrak{n}}^{1,F}\,(0,0,-1,0^n)^t = 0\\
&(\nabla^{\tilde{\omega}^{nor}}_{\xi}\nu - \nabla^{\tilde{\omega}}_{\xi}\nu)(u) \subseteq \tilde{\mathfrak{n}}^{1,F}\,(0,1,-1,0^n) \subseteq (0,0,0,\mathbb{R}^n) \quad \text{and}\\
&(\nabla^{\tilde{\omega}^{nor}}_{\tilde{V}}\nu - \nabla^{\tilde{\omega}}_{\tilde{V}}\nu)(u) = 0
\end{aligned}$$

for all $\xi \in \mathfrak{X}(\tilde{M})$ and $u \in \mathcal{G}$. Note that

$$\mathcal{G} \times_H (0,0,0,\mathbb{R}^n) = \mathrm{Ann}(\eta) \subseteq \tilde{F}^* \subseteq \tilde{\mathcal{T}}^*.$$

Combining all above relations, we conclude that

$$\nabla^{\tilde{\omega}^{nor}}\mu = 0, \quad \nabla^{\tilde{\omega}^{nor}}\nu \in \Omega^1(\tilde{M}, \mathrm{Ann}(\eta)) \quad \text{and} \quad \nabla^{\tilde{\omega}^{nor}}_{\tilde{V}}\nu = 0,$$

as desired.

Moreover, recall from Lemma 25 that $(j^*\tilde{\kappa}^{nor})(\tilde{V}, T\tilde{M}) \subseteq \mathcal{G} \times_H \tilde{\mathfrak{n}}^{1,F}$. Observe that $\mathcal{G} \times_H \tilde{\mathfrak{n}}^{1,F} = \mathcal{N}^{1,F}$. This completes the proof of the corollary. □

With the consideration of Proposition 3 and our knowledge about $\tilde{\partial}^*\tilde{\kappa}$ from the proof of Proposition 24, we may modify Theorem 17 carefully to obtain a characterization of all almost Grassmann structures of type $(2, n+1)$ that are locally the result of the Fefferman-type construction on an almost Grassmann structure of type $(2, n)$ via our Fefferman-type construction.

Now let $(\tilde{M}, \tilde{E}, \tilde{F})$ be an almost Grassmann structure of type $(2, n+1)$ and $(\tilde{\mathcal{G}} \to \tilde{M}, \tilde{\omega}^{nor})$ be the normal parabolic geometry of type $(\tilde{G}, \tilde{P})$ canonically associated to it. Denote the curvature of $\tilde{\omega}^{nor}$ by $\tilde{\kappa}^{nor}$. Recall the short exact sequences

$$0 \to \tilde{E} \to \tilde{\mathcal{T}} \to \tilde{F} \to 0 \quad \text{and} \quad 0 \to \tilde{F}^* \to \tilde{\mathcal{T}}^* \to \tilde{E}^* \to 0$$

where $\tilde{\mathcal{T}}$ is the standard tractor bundle of $(\tilde{\mathcal{G}} \to \tilde{M}, \tilde{\omega}^{nor})$ with dual $\tilde{\mathcal{T}}^*$. Each of the standard tractor bundle $\tilde{\mathcal{T}}$ and the standard cotractor bundle $\tilde{\mathcal{T}}^*$ are equipped with the tractor connection induced by $\tilde{\omega}^{nor}$. We denote them by $\nabla^{\tilde{\omega}^{nor}}$. Moreover, every nowhere vanishing section

$$\eta \in \Gamma(\tilde{F})$$

identifies $\tilde{E}^*$ with the distribution

$$\tilde{V} = \tilde{E}^* \otimes \eta \subseteq T\tilde{M}.$$

We also obtain a refined filtration

$$\tilde{E} \subseteq \mathrm{span}(\eta) + \tilde{E} \subseteq \tilde{\mathcal{T}}$$

of ranks $2 \leq 3 \leq n+3$. Its dual filtration is

$$\mathrm{Ann}(\eta) \subseteq \tilde{F}^* \subseteq \tilde{\mathcal{T}}^*$$

of ranks $n \leq n+1 \leq n+3$. Here, $\mathrm{Ann}(\eta)$ denotes the annihilators of $\eta$ in $\tilde{F}^*$. Let

$$\tilde{\mathcal{N}}^{1,F} = \{\Psi \in \mathfrak{sl}(\tilde{\mathcal{T}}) : \Psi(\mathrm{span}(\eta) + \tilde{E}) = 0 \text{ and } \Psi(\tilde{\mathcal{T}}) \subseteq \mathrm{span}(\eta) + \tilde{E}\}.$$

**Theorem 27** *Notation as above. Then $(\tilde{M}, \tilde{E}, \tilde{F})$ is locally the result of the Fefferman-type construction on an almost Grassmann structure of type $(2, n)$ if and only if the following hold.*

*(i)*
- *There is a parallel standard tractor $\mu \in \Gamma(\tilde{\mathcal{T}})$ with nowhere vanishing projection $\eta \in \Gamma(\tilde{F})$.*
- *There is a standard cotractor $\nu \in \Gamma(\tilde{\mathcal{T}}^*)$ such that*

$$\nabla^{\tilde{\omega}^{nor}} \nu \in \{\Phi \in \Omega^1(\tilde{M}, \tilde{F}^*) : \Phi(\tilde{V}) = 0 \textit{ and } \Phi(T\tilde{M}) \subseteq \mathrm{Ann}(\eta)\}.$$

*Moreover, there holds $\nu(\mu) = 1$ and the projection $\varphi \in \Gamma(\tilde{E}^*)$ of $\nu$ is nowhere vanishing.*

*(ii)*
- $\tilde{\kappa}^{nor}(\tilde{V}, \tilde{V}) = 0$ *and*
- $\tilde{\kappa}^{nor}(\tilde{V}, T\tilde{M}) \subseteq \tilde{\mathcal{N}}^{1,F}$.

***Proof*** That (i) and (ii) hold for a given Fefferman-type construction follows from Lemma 25 together with Corollary 26. It remains to show the converse direction.

Indeed, consider the two grading decompositions of $\tilde{\mathfrak{g}}$ given in Convention 6 and Convention 23 and the $H$-modules $\mathbb{E}$ and $\mathbb{F}$ given in the proof of Lemma 25. Let $(\tilde{M}, \tilde{E}, \tilde{F})$ be an almost Grassmann structure of type $(2, n+1)$, $(\tilde{\mathcal{G}} \to \tilde{M}, \tilde{\omega}^{nor})$ the normal parabolic geometry of type $(\tilde{G}, \tilde{P})$ canonically associated to it and $\mu$ a tractor and $\nu$ a cotractor of $(\tilde{\mathcal{G}} \to \tilde{M}, \tilde{\omega}^{nor})$ satisfying the properties in (i) and (ii) of the statement of the theorem. We show the following.

- There is a reduction $j : \mathcal{G} \to \tilde{\mathcal{G}}$ to structure group $H$ and a Cartan geometry

$$(\mathcal{G} \to \tilde{M}, \omega)$$

of type $(G, H)$ such that

$$i' \circ \omega - j^*\tilde{\omega}^{nor} = \phi \circ (j^*\tilde{\omega}^{nor})$$

for some $\phi \in \Gamma(\mathcal{G} \times_H \mathbb{E})$ where $\mathbb{E}$ is as given in the proof of Lemma 25. Moreover,
- $(\mathcal{G} \to \tilde{M}, \omega)$ is locally the result of the correspondence space modulo $L(\mathfrak{g}/\mathfrak{p}, \mathfrak{p}_+)$ of a normal parabolic geometry of type $(G, P)$; see Definition 2.

Note that two Cartan connections of type $(\tilde{G}, \tilde{P})$ on $\tilde{\mathcal{G}} \to \tilde{M}$ with difference in $\tilde{\mathfrak{p}}$ have the same underlying almost Grassmann structure; see, e.g., [11, Proposition 3.1.10 (1)]. In particular, since $\mathbb{E} \subseteq \tilde{\mathfrak{p}}_+ \otimes \tilde{\mathfrak{n}}^{1,F}$ and $i'(\mathfrak{p}_+) \subseteq \tilde{\mathfrak{n}}^{1,F} \subseteq \tilde{\mathfrak{p}}$, the above two assertions imply that $(\tilde{M}, \tilde{E}, \tilde{F})$ is locally isomorphic to the result of the Fefferman-type construction.

Indeed, note that a section $\mu$ of the standard tractor bundle and a section $\nu$ of the standard cotractor bundle such that $\nu(\mu) = 1$ and that their projections $\eta$ and $\varphi$, respectively, are nowhere vanishing uniquely determine a reduction

$$j : \mathcal{G} \to \tilde{\mathcal{G}}$$

to structure group $H$ such that

$$\mu = \mathcal{G} \times_H (0, 0, -1, 0^n)^t \quad \text{and} \quad \nu = \mathcal{G} \times_H (0, 1, -1, 0^n).$$

This is proven in the third paragraph in the proof of Theorem 17. In particular, there holds

$$\tilde{V} = \tilde{E}^* \otimes \eta = \mathcal{G} \times_H \frac{\tilde{\mathfrak{n}}^E_{-1} \oplus \tilde{\mathfrak{p}}}{\tilde{\mathfrak{p}}} \quad \text{and} \quad \mathcal{G} \times_H \tilde{\mathfrak{n}}^{1,F} = \tilde{\mathcal{N}}^{1,F}. \tag{41}$$

Moreover, the assumptions

$$\nabla^{\tilde{\omega}^{nor}} \mu = 0, \quad \nabla^{\tilde{\omega}^{nor}}_{\tilde{V}} \nu = 0 \quad \text{and} \quad \nabla^{\tilde{\omega}^{nor}} \nu \in \Omega^1(\tilde{M}, \mathrm{Ann}(\eta))$$

in the statement of the theorem imply that

$$j^* \tilde{\omega}^{nor} \in \Omega^1(\mathcal{G}, \tilde{\mathfrak{g}})$$

takes values in

$$\begin{aligned} &\{A \in \tilde{\mathfrak{g}} : A \cdot (0, 0, -1, 0^n)^t = 0 \text{ and } (0, 1, -1, 0^n) \cdot A \in (0, 0, 0, \mathbb{R}^n)\} \\ &\quad = i'(\mathfrak{g}) + \tilde{\mathfrak{n}}^{1,F} \end{aligned}$$

where $\tilde{\mathfrak{n}}^{1,F}$ is given in Convention 23. Moreover, the restriction of $j^* \tilde{\omega}^{nor}$ to the maximal distribution projecting to $\tilde{V} = \tilde{E}^* \otimes \eta \subseteq T\tilde{M}$ takes values in

$$\begin{aligned} &\{A \in \tilde{\mathfrak{g}} : A \cdot (0, 0, -1, 0^n)^t = 0 \text{ and } (0, 1, -1, 0^n) \cdot A = 0\} \\ &\quad = i'(\mathfrak{g}). \end{aligned}$$

Let

$$q : i'(\mathfrak{g}) + \tilde{\mathfrak{n}}^{1,F} \twoheadrightarrow \mathfrak{g}$$

be given by deleting the third row and the third column of $\tilde{\mathfrak{g}}$. Observe that $q|_{i'(\mathfrak{g})}$ and $q|_{\tilde{\mathfrak{n}}^{1,F}}$ are $H$-equivariant so that, in particular,

$$\omega = q \circ j^* \tilde{\omega}^{nor} \in \Omega^1(\mathcal{G}, \mathfrak{g})$$

is $H$-equivariant. Since $q \circ i'$ restricts to the identity map on $\mathfrak{h}$, $\omega$ reproduces the generator of fundamental vector fields. Finally, since $q$ is surjective, $\omega$ restricted to every tangent space is a linear isomorphism. Hence

$$(\mathcal{G} \to \tilde{M}, \omega) \tag{42}$$

is a Cartan geometry of type $(G, H)$. By construction,

$$\omega - j^* \tilde{\omega}^{nor} = \phi \circ (j^* \tilde{\omega}^{nor}) \in \Omega^1(\mathcal{G}, \mathfrak{g})^H$$

for some

$$\phi \in \Gamma(\mathcal{G} \times_H \mathbb{E}).$$

Denote by $\kappa$ the curvature of $\omega$. As mentioned in the proof of Lemma 25, a computation using the formula of Cartan curvature shows that

$$i' \circ \kappa - j^* \tilde{\kappa}^{nor} \in \Gamma(\mathcal{G} \times_H \mathbb{F}).$$

In particular, it follows from the assumptions in (ii) of the statement of the theorem that

$$\kappa(\tilde{V}, \tilde{V}) = 0 \quad \text{and}$$
$$\kappa(\tilde{V}, T\tilde{M}) \subseteq \mathcal{G} \times_H \mathfrak{p}_+.$$

Note that $\mathfrak{p}_+ = (i')^{-1}(\tilde{\mathfrak{n}}^{1,F})$. Hence $(\mathcal{G} \to \tilde{M}, \omega)$ is locally isomorphic to the correspondence space of some Cartan geometry of type $(G, P)$ modulo $L(\mathfrak{g}/\mathfrak{p}, \mathfrak{p}_+)$; see Proposition 3. It remains to show that such a Cartan geometry of type $(G, P)$ may be chosen to be normal. Without loss of generality, we may assume that there is a parabolic geometry

$$(\mathcal{G} \to M, \omega')$$

of type $(G, P)$ such that $\mathcal{G}/H = \tilde{M}$ and $\omega' = \omega$ modulo $L(\mathfrak{g}/\mathfrak{p}, \mathfrak{p}_+)$. Equivalently,

$$\omega' - \omega = \phi' \circ \omega \in \Omega^1(\mathcal{G}, \mathfrak{g})^H$$

for some

$$\phi' \in \Gamma(\mathcal{G} \times_H L(\mathfrak{g}/\mathfrak{p}, \mathfrak{p}_+)) \subseteq \Gamma(\mathcal{G} \times_H \mathbb{E}).$$

We have viewed $(\mathfrak{g}/\mathfrak{p})^* \subseteq (\mathfrak{g}/\mathfrak{h})^* \cong (\tilde{\mathfrak{g}}/\tilde{\mathfrak{p}})^*$. Note that $i'(\mathfrak{p}_+) \subseteq \tilde{\mathfrak{n}}^{1,F}$. Moreover, since $j^* \tilde{\omega}^{nor}$ and $\omega$ induce the same identification

$$T\tilde{M} \cong \mathcal{G} \times_H \mathfrak{g}/\mathfrak{h} = \mathcal{G} \times_H \tilde{\mathfrak{g}}/\tilde{\mathfrak{p}},$$

we conclude that

$$i' \circ \omega' - j^*\tilde{\omega}^{nor} = (\phi + \phi') \circ (j^*\tilde{\omega}^{nor})$$

and hence $i' \circ \kappa' - j^*\tilde{\kappa}^{nor} \in \Gamma(\mathcal{G} \times_H \mathbb{F})$. Since $j^*\tilde{\kappa}^{nor} \in \Gamma(\mathcal{G} \times_H \mathbb{F})$, there holds

$$i' \circ \kappa' \in \Gamma(\mathcal{G} \times_H \mathbb{F}).$$

Now, denote by $\partial^*$ and $\tilde{\partial}^*$ the Kostant codifferential of $(\mathfrak{g}, \mathfrak{p})$ and of $(\tilde{\mathfrak{g}}, \tilde{\mathfrak{p}})$, respectively. Let $\beta : \tilde{\mathfrak{g}} \to \mathfrak{g}$ be as given in (28) in the proof of Theorem 17, i.e., $\beta$ is defined by deleting the third row and summing up the second and the third column of $\tilde{\mathfrak{g}}$. Since $\partial^*$ is the restriction of the Kostant codifferential on $(\mathfrak{g}, \mathfrak{q})$, the relation (28) between $\tilde{\partial}^*$ and $\partial^*$ still holds. Hence

$$\partial^*\kappa' = \beta \circ \tilde{\partial}^*(i' \circ \kappa').$$

Since

$$\tilde{\partial}^*(i' \circ \kappa') \in \tilde{\partial}^*\,\mathbb{F} \subseteq \mathbb{E}$$

and the composition of an element in $\mathbb{E}$ with $\beta$ lies in $L(\mathfrak{g}/\mathfrak{p}, \mathfrak{p}_+)$, we conclude that

$$\partial^*\kappa' \in \Omega^1_{hor}(\mathcal{G}, \mathfrak{p}_+)^P.$$

Therefore, there exists a normal parabolic geometry

$$(\mathcal{G} \to M, \omega'_{nor})$$

of type $(G, P)$ such that $\omega'_{nor} = \omega'$ modulo $L(\mathfrak{g}/\mathfrak{p}, \mathfrak{p}_+)$. This completes the proof of the theorem. □

***Remark 28*** Notation as in Theorem 27. Note that by Bernstein-Gelfand-Gelfand (BGG) theory, $\nabla^{\tilde{\omega}^{nor}} \nu \in \Omega^1(\tilde{M}, \tilde{F}^*)$ means that $\nu$ is the image of $\varphi$ under the splitting operator associated to $\tilde{\mathcal{T}}^*$ and $\varphi$ is a solution of the so-called first BGG operator of $\tilde{\mathcal{T}}^*$. Similarly, $\nabla^{\tilde{\omega}^{nor}} \mu = 0$ means that $\mu$ is the image of $\eta$ under the splitting operator associated to $\tilde{\mathcal{T}}$ and $\eta$ is a normal solution of the first BGG operator of $\tilde{\mathcal{T}}$. In fact, analogous to Corollary 19, the conditions in Theorem 27 (i) can be expressed as conditions on sections $\varphi \in \Gamma(\tilde{E}^*)$ and $\eta \in \Gamma(\tilde{F})$ in terms of a Weyl connection. To do so, one applies [18, Proposition 5 and Lemma 7] in the situation of Theorem 27 (i).

**Funding** Open access funding provided by University of Vienna.

**Data Availability** Data sharing is not applicable to this article as no datasets were generated or analysed during the current study.

**Publisher's Note** Springer Nature remains neutral with regard to jurisdictional claims in published maps and institutional affiliations.